\documentclass[10pt]{article}
\usepackage{pdfsync, caption} 
\usepackage{listings}
\usepackage[latin1]{inputenc}
\usepackage{longtable,tabu} 
\usepackage{color} %red, green, blue, yellow, cyan, magenta, black, white
\definecolor{mygreen}{RGB}{28,172,0} % color values Red, Green, Blue
\definecolor{mylilas}{RGB}{170,55,241}
\usepackage[pagebackref=true, colorlinks=true, linkcolor=blue, citecolor=blue, linkcolor=blue, anchorcolor=red, urlcolor=blue]{hyperref}
\usepackage{url, amsmath,enumerate,fancyhdr,amssymb, amsthm, 
pstricks, epsf, lscape, ifpdf, graphicx, float}
\usepackage{algorithm}
\usepackage{color}
\usepackage{algpseudocode} 
\usepackage{tabularx, makecell, multirow} 
\usepackage{lscape}
\usepackage[numbers,sort]{natbib}

\algnewcommand\algorithmicinput{\textbf{Input:}}
\algnewcommand\INPUT{\item[\algorithmicinput]}
\algnewcommand\algorithmicinitialization{\textbf{Initialization:}}
\algnewcommand\INITIALIZATION{\item[\algorithmicinitialization]}

\usepackage[margin=3cm]{geometry}

\renewcommand{\arraystretch}{1.3}

\mathchardef\mhyphen="2D

\newtheorem{Definition}{Definition}

\newtheorem{Example}{Example}
\newtheorem{Proposition}{Proposition}
\newtheorem{Lemma}{Lemma}
\newtheorem{Theorem}{Theorem}
\newtheorem{Corollary}{Corollary}
\newtheorem{Remark}{Remark}
\newtheorem{Assumption}{Assumption}

\newcommand{\A}{\mathcal A}

\newcommand{\err}{\operatorname{err}}

\newcommand{\eps}{\epsilon} 
\newcommand{\bpx}{\begin{pmatrix}}
\newcommand{\epx}{\end{pmatrix}}
\newcommand{\bbx}{\begin{bmatrix}}
\newcommand{\ebx}{\end{bmatrix}}

\newcommand{\bdef}{\begin{Definition}} 
\newcommand{\comment}[1]{}
\newcommand{\co}[1]{}

\newcommand{\coab}[1]{}

\newcommand{\norm}[1]{\parallel \!\! #1 \!\! \parallel}
\newcommand{\sym}[1]{{\cal S}^{#1}}
\newcommand{\psd}[1]{{\cal S}_+^{#1}}
\newcommand{\rad}[1]{\mathbb{R}^{#1}}

\newcommand{\beq}{\begin{equation}}
\newcommand{\eeq}{\end{equation}}
\newcommand{\beqa}{\begin{eqnarray}}
\newcommand{\eeqa}{\end{eqnarray}}
\newcommand{\ba}{\begin{array}}
\newcommand{\ena}{\end{array}}
\newcommand{\bac}{\begin{array}{ccccccccccc}}
\newcommand{\eac}{\end{array}}
\newcommand{\bprop}{\begin{Proposition}}
\newcommand{\eprop}{\end{Proposition}}

\newcommand{\beqast}{\begin{eqnarray*}}
\newcommand{\eeqast}{\end{eqnarray*}}
\newcommand{\benum}{\begin{enumerate}}
\newcommand{\eenum}{\end{enumerate}}
\newcommand{\bit}{\begin{itemize}}
\newcommand{\eit}{\end{itemize}}
\newcommand{\bth}{\begin{Theorem}}
\newcommand{\enth}{\end{Theorem}}
\newcommand{\ble}{\begin{Lemma}}
\newcommand{\ele}{\end{Lemma}}
\newcommand{\bex}{\begin{Example}}
\newcommand{\eex}{\end{Example}}
\newcommand{\bcor}{\begin{Corollary}}
\newcommand{\ecor}{\end{Corollary}}
\newcommand{\brem}{\begin{Remark}}
\newcommand{\erem}{\end{Remark}}
\newcommand{\bass}{\begin{Assumption}}
\newcommand{\eass}{\end{Assumption}}

\renewcommand{\arraystretch}{1.2}

\newcommand{\val}{\operatorname{val}}

\newcommand{\bsmx}{\begin{small} \begin{pmatrix}}
\newcommand{\esmx}{\end{pmatrix} \end{small}}

\title{\Large Sieve-SDP:  a simple facial reduction algorithm to preprocess semidefinite programs} 
\author{Yuzixuan Zhu$^{\dagger}$ \hspace{1cm} G\'{a}bor  Pataki$^{\dagger}$\thanks{Corresponding author}\hspace{1cm} Quoc Tran-Dinh\thanks{Y. Zhu, G. Pataki,  and Q. Tran-Dinh are with the Department of Statistics and Operations Research, University of North Carolina at Chapel Hill.\newline
\textit{Address:} Hanes Hall, Chapel Hill, NC 27599-3260. \textit{Email:} zyzx@live.unc.edu, gabor@unc.edu, quoctd@email.unc.edu. 
}
}

\begin{document}

\maketitle

\begin{abstract}
We introduce Sieve-SDP, a simple facial reduction algorithm to preprocess semidefinite programs (SDPs). Sieve-SDP inspects the constraints of the problem to detect lack of strict feasibility, deletes redundant rows and columns, and reduces the size of the variable matrix. It often detects infeasibility.  
It does not rely on any optimization solver: the only subroutine it needs is Cholesky factorization, 
hence it can be implemented in a few lines of code  in machine precision. 
We present extensive computational results on several problem collections  from the literature, with many SDPs coming from polynomial optimization. 
\end{abstract}



%%%%%%%
% \section{Introduction and the preprocessing algorithm}
%%%%%%%

\section{Introduction and the preprocessing algorithm}\label{sec:intro}

Consider  a semidefinite programming problem (SDP) in the form 
%Semidefinite programs (SDPs) generalize linear programs
\begin{equation} \label{p} \tag{\mbox{$P$}} 
\left.\begin{array}{rrcl}
\displaystyle\inf_{X} & C \bullet X \\
%s.t.  &                 \\
\mathrm{s.t.} ~& A_{i}\bullet X & = & b_{i} \,\,(i=1,\ldots,m), \\
& X & \succeq & 0,
\end{array}\right.
\end{equation}
where the $A_i$ and $C$ are $n \times n$ 
symmetric matrices, the $b_i$ are scalars, $X \succeq 0$ means that $X$ is 
in $\psd{n}, \,$ the set of symmetric, 
positive semidefinite (psd) matrices,  
and the $\bullet$ inner product of symmetric matrices is 
the trace of their regular product.

SDPs are some of the most versatile, useful, and widespread optimization problems of the last three decades. They find  applications in 
control theory, integer programming,   and combinatorial optimization, to name just a few areas. 
Several good  solvers are available to solve SDPs (see for example 
\cite{TutTohTod:03, Sturm:99, fujisawa2002sdpa, fujisawa2008sdpaGMP, zhao2010newton,  KocvaraStingl:2003, BurerMonteiro:2003, BurerMonteiroZhang:2002, Mosek:14});  among these, Mosek \cite{Mosek:14}  is commercially available.

SDPs -- as all optimization problems -- often have  redundant variables and/or constraints. 
The redundancy we address is lack of {\em strict feasibility}, 
%when \eqref{p} is not strictly feasible, 
i.e., when there is no feasible positive definite $X$ in \eqref{p}. 
When \eqref{p} is not strictly feasible, % strong duality may fail, 
%i.e.,
 the optimal value of \eqref{p} and of its dual may differ, and the latter may not be attained\footnote{More precisely, when \eqref{p} {\em is} strictly feasible, {\em strong duality} holds between \eqref{p} and its dual, i.e., their  values agree and the latter is attained.}. 
 Hence,  when attempting to solve such an SDP, solvers often struggle, or fail.  

It is, of course, useful to detect lack of strict feasibility 
in a preprocessing stage.      
This paper describes a very simple preprocessing algorithm for SDPs, called Sieve-SDP, which belongs to the class of facial reduction algorithms  \cite{BorWolk:81, WakiMura:12, Pataki:00B, Tuncel:11, Pataki:13, KrisWolk:10, Dima:15, drusvyatskiy2014noisy, PerPar:14}. 
%\cite{drusvyatskiy2014noisy} \cite{Dima:15} \cite{KrisWolk:10} 
 Sieve-SDP  
can 
detect lack of strict feasibility, reduce the size of the problem, and can 
be implemented in a few lines of code in machine precision.

To motivate our algorithm, let us  consider  an example:
\bex \label{ex1} 
The SDP instance $($with an arbitrary objective function$)$
\begin{equation}  \label{mot-ex}
\left.\begin{array}{rcl} \vspace{.2cm} 
\begin{pmatrix} 1 & 0 & 0 \\
0 & 0 & 0 \\
0 & 0 & 0 
\end{pmatrix} \, \bullet  X & = & 0 \\  
\begin{pmatrix} 0 & 0 & 1 \\
0 & 1 & 0 \\
1 & 0 & 0 
\end{pmatrix} \, \bullet \, X & = & -1 \\
 X & \succeq & 0,
\end{array}\right.
\end{equation}
is infeasible. Indeed, suppose $X = (x_{ij})_{i,j=1}^3$ is feasible in (\ref{mot-ex}). 
Then $x_{11}=0,$   hence the first row and column of $X$ are zero
by positive semidefiniteness, so  the second constraint implies $x_{22} = -1, \,$ which is a contradiction.
\eex
Note that if we replace $-1$ in the second constraint of 
(\ref{mot-ex}) by a positive number, then (\ref{mot-ex}) can be restated 
over the set of psd matrices with first row and column equal to zero. Thus, even if we do not detect infeasibility, such 
preprocessing is still useful.

\begin{figure}[H]
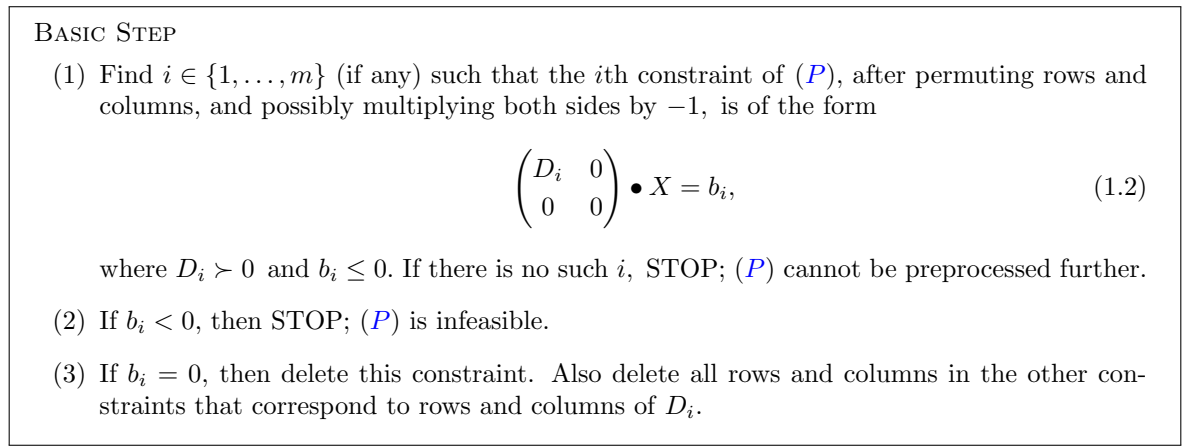

\vspace{-0.5ex}
\begin{center}
	\framebox[\textwidth]{\parbox{0.95\textwidth}{ 
			{\sc Basic Step}
			\begin{enumerate}
				\vspace{-1ex}
				\item \label{findD} Find  $i \in \{1, \dots, m \}$ (if any) such that the $i$th constraint of \eqref{p},
				after permuting rows and columns, and possibly multiplying both sides by $-1, \,$ 
				is of the form  
				\beq \label{basic} 
				\bpx D_i & 0 \\ 0 & 0 \epx \bullet X = b_i, 
 				\eeq
				where $D_i \succ 0  \, $ and $b_i \leq 0.$ If there is no such $i, \,$ STOP; \eqref{p} cannot be preprocessed  further. 
				\item \label{checkbi<0} If $b_i < 0$, then STOP; \eqref{p} is infeasible.
				\item \label{checkbi=0} If $b_i = 0$, then delete this constraint. Also delete 
				all rows and columns in the other constraints that correspond to rows and columns of $D_i.$ 
				\vspace{-1ex}
			\end{enumerate}
		}
	}
\vspace{-1ex}	
\caption{The Basic Step  of Sieve-SDP}\label{figure-basicstep} 
\end{center}
\vspace{-1ex}
\end{figure}

Our algorithm Sieve-SDP repeats the Basic Step shown in Figure \ref{figure-basicstep}. 
Hereafter  $D \succ 0$ means  that a symmetric matrix  $D$ is positive definite.
		
		\bex (Example \ref{ex1} continued)  When we first execute the Basic Step on \eqref{mot-ex}, 
	we find the first constraint, delete it, and also delete the first row and column from the second constraint matrix. Next, we find the constraint
	$$
	\bpx 1 & 0 \\ 0 & 0 
	\epx \bullet X = -1,
	$$
	and declare  that (\ref{mot-ex}) is infeasible. 
	\eex
		
\begin{figure}[!htp]
	\centering
	\includegraphics[width = 6cm]{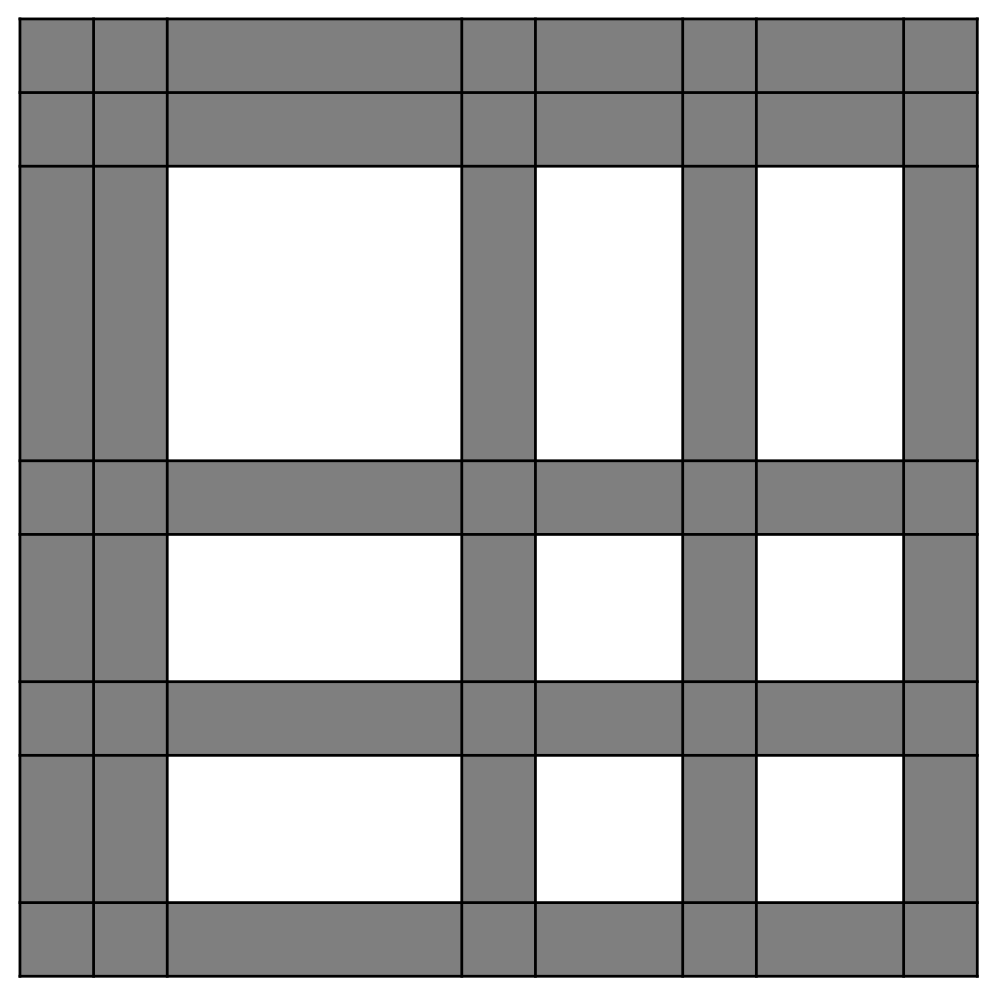}
\caption{The sieve structure}
\label{figure:sieve} 
\end{figure}

We call our  algorithm Sieve-SDP, since by shading the deleted rows and columns in the 
variable matrix $X$ (and the $A_i$) we obtain a sieve-like structure: 
see Figure \ref{figure:sieve}.
		
Sieve-SDP is easy to implement and fast: it only needs   
an incomplete Cholesky factorization subroutine 
to check positive definiteness, and 
we can delete rows and columns using fast matrix operations. Even the worst case complexity of Sieve-SDP is reasonable: an easy calculation shows that it can fully preprocess \eqref{p} using 
$O(\min \{ m, n \}  n^3 m)$ arithmetic operations.

Sieve-SDP is a heuristic: it does not always detect infeasibility, or lack of strict feasibility.
For example, it will not work on problem \eqref{mot-ex}, if we  apply a similarity transformation 
 $T^{\top}(\cdot)T$ to all $A_i, \,$ where $T$ is a random invertible matrix.
 
Sieve-SDP is very simple, and  easy to ``fool", thus, 
%Given  its simplicity, and how easily it is ``fooled",  
it is natural to ask whether it works in practice. 
So the main research question we address,  and answer in the affirmative,  is:
\begin{itemize}
	\item Can Sieve-SDP  help us compute more accurate solutions 
and reduce the computing time on a broad range of SDPs?
\end{itemize}

\paragraph{\bf Related work:} Sieve-SDP belongs to the family of {\em facial reduction algorithms}, which we now describe. 
When \eqref{p} is not strictly feasible, one can replace the constraint $X \in \psd{n}$ by 
$$
X \in F,
$$
where $F$ is a proper face of $\psd{n}. \footnote{That is, $F \neq \psd{n}, \,$ $F$ is convex, and $X, Y \in \psd{n},  \, \frac{1}{2}(X+Y) \in F$ implies that $X$ and 
	$Y$ are in $F.$}$ 
Since any such face can be written as (see e.g. \cite{Pataki:00A}) 
\begin{equation} \label{VSVT} 
F = V\psd{r}V^{\top}, 
\end{equation}
where $r < n \,$ and $V$ is an $n \times r$ matrix, the reduced problem can be restated over a smaller semidefinite cone. Facial reduction algorithms
-- for more general conic  programs -- originated in the papers \cite{BorWolk:81, BorWolk:81B}. Later simplified, more easily implementable variants were  given in \cite{Pataki:00B, WakiMura:12, Pataki:13}, and in \cite{Tuncel:11} for the SDP case.  
A recent, very concise version with a short proof of convergence is in \cite{LiuPataki:17}.

Facial reduction algorithms, when applied to \eqref{p}, find  the face $F$ by solving a sequence of SDP
%When we apply a facial reduction algorithm to \eqref{p}, it finds the face $F$ by solving a sequence of SDP
 subproblems, which may be as hard to 
solve as \eqref{p} itself. Thus one is led  to  seek simpler alternatives.

Simplified and implementable versions of  facial reduction are described in 
\cite{PerPar:14}. The algorithms in \cite{PerPar:14} reduce the feasible set of (\ref{p}) (or of an SDP in a different shape) 
by  solving  linear programs instead of SDPs. 
Thus they do not find {\em all} reductions, but still simplify the SDPs in many cases.
They are available as  public domain codes, and 
we will compare them  with Sieve-SDP in Section \ref{section:implement}. A facial reduction algorithm embedded in an interior point method was implemented in \cite{permenter2017solving}. 

We next review facial reduction algorithms that work by simply inspecting 
constraints. For example,  \cite{Friberg:16} notes that
if
$$A \bullet X = 0$$ is a constraint in 
\eqref{p} with $A \succeq 0,  \,$ then we can restrict $X$ to belong to a face of the form 
(\ref{VSVT}) (where $V$ spans the nullspace of $A$). 
A similar idea was used in \cite{KrisWolk:10} to reduce Euclidean Distance Matrix completion problems. For a rigorous derivation of the algorithm in 
\cite{KrisWolk:10} see \cite{Dima:15}, which used an intermediate step of analyzing 
the semidefinite completion problem. 
For followup work, see \cite{ drusvyatskiy2014noisy} on the noisy version of the same problem, and \cite{tanigawa2017singularity} for a more theoretical study.

We finally mention two very accurate  SDP solvers, which do not rely on facial reduction. The first is SDPA-GMP \cite{fujisawa2008sdpaGMP}, which  %uses the GMP library  and 
 computes solutions of \eqref{p} and of its dual using several hundred digits of accuracy. We will use SDPA-GMP in later sections
 to check 
 the accuracy of the solutions computed by Sieve-SDP and Mosek. 
The SPECTRA solver \cite{HenrionNaldi:2016}   computes a  feasible solution of \eqref{p} (if one exists) in exact arithmetic. Although these solvers cannot handle large SDPs, they can solve small ones very accurately.

Sieve-SDP differs in several aspects from previously proposed facial reduction 
algorithms:
\begin{itemize}
	\item  	%It needs only  Cholesky factorization as a subroutine 	
	%and 
	%it does not rely on any optimization solver, like the algorithms in \cite{PerPar:14}.
	
	It needs only  Cholesky factorization as a subroutine
	and, unlike the algorithms in \cite{PerPar:14}, 	it does not rely on any optimization solver. 
	
%	It needs only  Cholesky factorization as a subroutine. Unlike the algorithms in \cite{PerPar:14}, 	
%	Sieve-SDP does not rely on any optimization solver. 
	%, like the algorithms in \cite{PerPar:14}.
\item It detects very simple redundancies, which are easy to explain 
even to a user not trained in 
optimization, and can help  him/her to better formulate other problems. 
\item As soon as Sieve-SDP  finds  a reducing constraint, 
it deletes this  constraint, and it also deletes redundant rows and columns from the other constraint matrices.
Hence errors do not accumulate. Thus Sieve-SDP is as accurate as   Cholesky factorization, which works  in machine precision  \cite[Theorem 23.2]{trefethen1997numerical}.
\item Sieve-SDP can also detect infeasibility.
\item It  is easy to run in a {\em safe mode} (explained in the next section) to even better 
safeguard %even better 
against numerical errors.
\item Finally, we present extensive computational results on general SDPs, 
which,  as  far as we know, are not yet  available for such a simple algorithm. 
\end{itemize}

The rest of the paper is organized as follows. In Section \ref{section:implement} 
we describe how we  implemented Sieve-SDP, the computational setup, and the 
criteria for comparison with competing codes. 
In this section we also give  a small  SDP with a positive duality gap (in Example \ref{ex2}), 
and show 
how to construct a pair of primal-dual solutions with arbitrarily small constraint violation and arbitrarily small duality gap. This example shows  that a solution with a {\em smaller} DIMACS error
(see \cite{Mittelmann:2003}) may be actually  {\em less accurate}. We also show that such a less accurate solution is actually computed by Mosek, one of the leading SDP solvers.

In Section \ref{section:detailed} we comment in detail on the results on some of the problems, and on the strengths and weaknesses of the preprocessors. For example, 
we examine whether they help to  find the correct solution of numerically difficult SDPs; and how fast they are on large scale problems. 

In Section \ref{section:summary} we summarize the preprocessing results, and 
 conclude the paper.
 
We have four  appendices. In Appendix  \ref{section:verydetailed} we 
  present very detailed 
 computational results on all  problems. 
 In Appendix \ref{section:matlabcode} we give the core 
 Matlab code of Sieve-SDP, containing  only about $65$ lines. 
 In Appendix \ref{section:dimacs} we provide  the definition of the DIMACS errors for completeness. 
 In Appendix 
 \ref{section:dual_recovery} we discuss the issue of recovering  an optimal solution of the dual of \eqref{p} from the optimal solution of the  dual of the reduced problem.

%%%%%%%
% \section{Implementation, setup for computational testing, and the issue of positive duality gaps}
%%%%%%%

\section{Implementation, setup for computational testing, codes used for comparison,  and the issue of positive duality gaps}\label{section:implement}

\subsection{Implementation and computing environment}\label{subsect:implementation}

We implemented our algorithm in Matlab R2015a, using the standard  
Cholesky factorization 
 (subroutine \texttt{chol}) 
to check positive definiteness. 

We  ran both Sieve-SDP and the competing preprocessors 
(which we describe in Subsection \ref{subsection-comparison})     
on a MacBook Pro with processor Intel Core i5 running at 
2.7GHz, and  8GB of RAM.

\subsection{Safe mode}\label{subsect:safemode}

To safeguard against numerical errors we use a {\em safe mode.} We set 
$$
\epsilon \, := 2^{-52} \, \approx \, 2.2204 \cdot 10^{-16} \, = \, {\rm the \, machine \, precision \, in \, Matlab.} 
$$   
In the Basic Step in Figure \ref{figure-basicstep}, 
if we find a constraint of type (\ref{basic}), then, 
 instead of checking $b_i < 0$ we check
whether 
$$
b_i < - \sqrt{\epsilon} \max \{ \, \norm{b}_{\infty}, 1 \} \, {\rm holds.} 
$$
If this test fails, then instead of checking $b_i = 0$ we check whether 
$$
b_i > - \epsilon \max \{ \, \norm{b}_{\infty}, 1 \} \, {\rm holds.} 
$$
Note that this step is correct, because  in the Basic Step we already ensured  $b_i \leq 0.$

\subsection{Preprocessors used for comparison}\label{subsection-comparison}

We compare Sieve-SDP with the algorithms proposed by Permenter and Parrilo in \cite{PerPar:14}. Their algorithms solve linear programming subproblems to reduce the size of an SDP. They can work either on the problem \eqref{p}, which we call the {\em primal}; or on its dual: 
\begin{equation} \label{d} \tag{\mbox{$D$}}   
\left.\begin{array}{rl}
\displaystyle\sup_y & \displaystyle\sum_{i=1}^m y_i b_i  \\
\mathrm{s.t.}   & \displaystyle\sum_{i=1}^m y_i A_i \preceq C.
\end{array}\right.
\end{equation}
They can use either diagonal, or diagonally dominant reductions (for details, see \cite{PerPar:14}).

Thus, there are four algorithms from \cite{PerPar:14} that we tested: 
pd1, pd2, dd1,  and dd2. 
Here pd1 stands for primal diagonal; pd2 for primal diagonally dominant;
dd1 for dual diagonal; and dd2 for dual diagonally dominant.

\brem{\rm  \label{remark:PP-primal-vs-dual} 
In the theoretical description of the algorithms in 
 \cite{PerPar:14} the SDP which  is called the primal is 
 actually  our dual (\ref{d}). %. their primal. 
 However, in  their implementation and their code posted on the github website,  
 their primal is the same as our primal  \eqref{p}.
}
\erem

\subsection{The datasets}\label{subsect:thedatasets}

We tested  Sieve-SDP and competing methods on five datasets, which contain 771 problems overall.

\begin{itemize} 
	\item The first is the dataset from \cite{PerPar:14}, which we call 
	the Permenter-Parrilo or PP dataset. 
	This dataset has 68 problems, whose original sources are 
	\cite{Waki:12, CheWolkSchurr:12, baston1969extreme, fawzi2016self, DIMACS,   diananda1962non, quist1998copositive, posa2013lyapunov,   boyd2014performance, waki2012strange, burton2014real, wagner2009criterion}.  Although a few problems in this dataset  are  randomly generated, 
	most come from applications.
\end{itemize}

The PP  dataset  contains SDPs that are notoriously difficult for 
solvers, and some  are known to be not strictly feasible. 
Hence we added the following  four  datasets to make our testing more comprehensive:

\begin{itemize}
	\item A dataset we obtained  from  Hans Mittelmann's website,  
	which we call the Mittelmann dataset.  
		This dataset contains 31 problems. 
	\item A collection of SDP relaxations of polynomial optimization problems based on the paper of 
		Dressler, Illiman, and de Wolff \cite{dressler2018approach}, which we call the 
		Dressler-Illiman-de Wolff dataset, or  DIW dataset for short. 
		This dataset has 155 problems.  
		\item A problem set kindly provided to us by Didier Henrion and Kim-Chuan Toh, which we call the Henrion-Toh dataset. This dataset contains 98 problems. 
	\item A problem set kindly provided to us by Kim-Chuan Toh,  whose 
		description is in \cite{sun2015convergent} and \cite{yang2015sdpnal}.  
		We call this dataset the Toh-Sun-Yang dataset, 
		and it has 419 problems.
\end{itemize}

From the PP dataset we excluded only two  problems: 
$\mathrm{copos\_5}$ and $\mathrm{cprank\_3}, \,$ since 
they were too large to be solved by Mosek on our computer.

Our datasets contain many different types of SDPs and, not surprisingly, 
the performance of the preprocessors on them 
varies widely. Many of our SDPs may be strictly feasible, 
and such SDPs 
could not be reduced by even  more sophisticated preprocessors.  
For example, in  the Toh-Sun-Yang dataset no  problems were reduced by the preprocessors.
Although this is a bit disappointing, Sieve-SDP and pd1  delivered the ``no reduction found" result very quickly, 
so it did  not hurt to preprocess.

%and on these SDPs 
%even more sophisticated preprocessors would not find reductions.  
%For example, on  the Toh-Sun-Yang dataset the preprocessors found no  reductions. 
%However, Sieve-SDP and pd1  took  only a negligible amount of time to deliver the ``no reduction %found" 
%result, so it did  not hurt to preprocess.

%and such SDPs 
%could not be reduced by even  more sophisticated preprocessors.  
%For example, in  the Toh-Sun-Yang dataset none of the problems were reduced. % by the preprocessors.
%Although this is a bit disappointing, Sieve-SDP and pd1 only took  a negligible amount of time to deliver the ``no reduction found" result, %so it did  not hurt to preprocess.

Yet, even in the datasets other than the PP dataset many SDPs {\em were}  reduced 
by some preprocessor.  In the Henrion-Toh dataset, pd1, pd2, and Sieve-SDP all 
reduced 18 problems, whereas  dd1 and dd2 reduced none. In the  Mittelmann dataset, 
pd1, pd2, and Sieve-SDP reduced 8 problems; dd1 and dd2 reduced none.

%Yet, even in the datasets other than the PP dataset 
%many SDPs {\em were}  reduced 
%by some preprocessor.  In the Henrion-Toh dataset 18 problems were reduced by  pd1, pd2, and Sieve-SDP, and none by dd1 and dd2. In %the  Mittelmann dataset  8 problems were reduced by 
%pd1, pd2, and Sieve-SDP, and none by dd1 and dd2. % reduced none.

Strikingly, in the DIW dataset 
Sieve-SDP  proved infeasibility of 59 problems out of 155, and reduced total solving time by a factor of more than a hundred! Pd1 did only slightly worse. % results.

\begin{figure}[!htbp]
	\begin{center}
		\includegraphics[scale = 0.4]{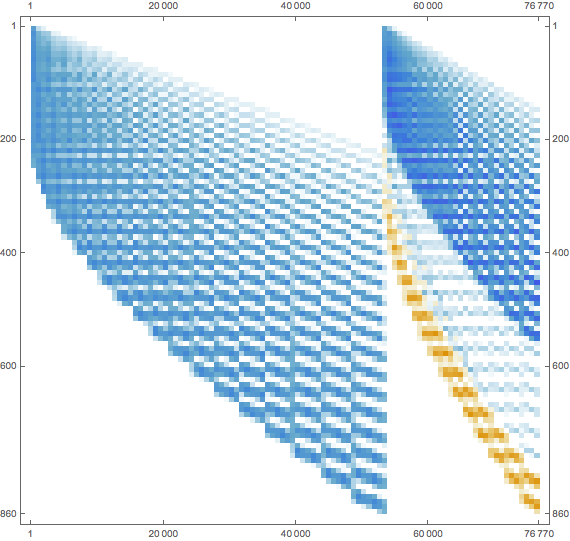} \hspace{0.2cm}\includegraphics[scale = 0.4]{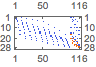} 
		\caption{Problem ``ex4.2\_order20'': size and sparsity before and after Sieve-SDP}\label{figure:DIW} 
	\end{center}
\vspace{-4ex}	
\end{figure}

We illustrate this point with Figure \ref{figure:DIW}, which shows the size and sparsity structure of the problem 
``ex4.2\_order20''\footnote{This SDP is from the DIW dataset.} 
before (on the left) and after (on the right) applying Sieve-SDP.
Each row in the displayed matrices corresponds to an $A_i$ matrix stretched out as a vector. Red dots correspond to  positive entries,  blue dots correspond to negative entries, and white areas correspond to zero entries.
%Sieve-SDP reduced this probkem from having $860$ rows to having only $28, \,$ then proved it to be infeasible. %at which point This SDP was found infeasible after 

\subsection{Internal format and input/output  format}\label{subsection-format}

Internally we store the $A_i$ matrices as an $n \times (nm)$ sparse matrix of the form
$$
\begin{pmatrix} A_1, A_2, \dots, A_m \end{pmatrix} 
$$
(i.e., the $A_i$ are stored side-by-side), and $C$ as an $n \times n$ sparse matrix. The input- and the output format of the preprocessors
is the widely used Mosekopt format.

\subsection{The choice of the SDP solver and LP solver}\label{subsect:computational-setup}

For all preprocessors we use %Both for Sieve-SDP and for the algorithms of \cite{PerPar:14} we use
Mosek 8.1.0.27  (from now on, simply ``Mosek") as SDP solver: we solve the SDPs with Mosek before and after preprocessing. 
We also solve the  linear programming (LP) subproblems in the algorithms of \cite{PerPar:14} by Mosek.
%We also use  Mosek as a linear programming (LP) solver (to solve the LP subproblems) in the algorithms of \cite{PerPar:14}.
We consider Mosek  as the best choice, since it is a reliable commercial SDP and LP solver, and it is being actively developed and improved.

Our settings are different from the ones used in \cite{PerPar:14}, where %in which 
 Sedumi \cite{Sturm:99} format is used as  input format, Mosek as  LP solver, and  Sedumi
 as  SDP solver.
With our settings  the algorithms of \cite{PerPar:14} work faster, because Mosek is much faster  than Sedumi. Although 
we must convert the data from Mosekopt   format to 
	Sedumi  format (to do the preprocessing), and then back (to solve the preprocessed problem with Mosek), 
the total conversion time is negligible: for each of pd1, pd2, dd1 and dd2 it is less than 100 seconds on {\em all}  771 SDPs. 
To be fair, in the detailed comparison tables of Appendix \ref{section:verydetailed} we list 
conversion time, and  preprocessing time separately. 

%Still,  in the detailed comparison tables of Appendix \ref{section:verydetailed} we list separately the time spent on 
 %conversion, and  the time spent   on preprocessing, to be fair. 

\subsection{Criteria for comparison}\label{subsection-criteria-comparison}

Let us recall the main question that we address in this paper: 
\begin{itemize}
\item 
Can Sieve-SDP  help us  compute 
more accurate solutions 
and reduce the computing time on a broad range of SDPs?
\end{itemize}

Thus,  we compare the preprocessors  based on the following three criteria: %are as follows, in order of priority:

\vspace{.1cm} 
\begin{enumerate} 
	  	\item  Do they  help detect infeasibility? If not, do they  help to  find a correct optimal solution?
	  	
	  	Precisely,  suppose that Mosek reports an incorrect optimal value of an SDP  before preprocessing. Does Mosek 
	  	find a correct optimal value after preprocessing ? (We assume that the optimal value of the SDP is known 
	mathematically.) 
	
	\item  Does preprocessing reduce computing time? 

This criterion is secondary, since 
preprocessing is often essential to compute any accurate solution: see Subsections 
\ref{subsect:compact} through \ref{subsect:example}. Thus,  we believe that we should always preprocess SDPs, 
as long as we can do this with very high precision, even if preprocessing 
 {\em increases} the solution time.
			\item 
			Does preprocessing improve numerical accuracy measured by the six DIMACS 
			errors \cite{Mittelmann:2003}\footnote{The description of the DIMACS errors  is given in Appendix \ref{section:dimacs}.}?  
		
		Let  
		$$\text{DIMACS}_\text{before} \,\, {\rm and} \,\, \text{DIMACS}_\text{after}$$ 
		be the largest absolute value of the  DIMACS errors before and after preprocessing, respectively. We say that a method 
		{\em improves}  the DIMACS error if it does not detect infeasibility and  
		$$\text{DIMACS}_\text{before} > {10}^{-6} ~~~~ {\rm and} ~~~~~ \dfrac{\text{DIMACS}_\text{after}}{\text{DIMACS}_\text{before}} < \dfrac{1}{10}.$$
		
		This last criterion must be taken with a grain of salt.
		While the DIMACS errors are very natural (they measure constraint violation and duality gap), 
		Example \ref{ex2} below shows that they do not always  
		measure accurately how good a solution is. In fact, a {\em larger}  
		DIMACS error may correspond to a {\em better} solution!
\end{enumerate}

\bex \label{ex2} Consider the SDP   
\beq \label{posgap} 
\left.\begin{array}{rrcl}  \vspace{.2cm}  
	\displaystyle\inf_X & \bpx 1 & 0 & 0 \\ 0 & 1 & 0 \\ 0 & 0 & 0 \epx \bullet X  \\ \vspace{.2cm}  
	\mathrm{s.t.} & \bpx 1 & 0 & 0 \\ 0 & 0 & 0 \\ 0 & 0 & 0 \epx \bullet X  & = & 0 \\ \vspace{.2cm}   
	& \bpx 0 & 0 & 1 \\ 0 & 1 & 0 \\ 1 & 0 & 0 \epx \bullet X  & = & 1 \\
	& X = (x_{ij}) & \succeq & 0,
\end{array}\right.
\eeq
and  its dual 
\beq \label{posgap-dual} 
\left.\ba{rl}
\displaystyle\sup_y  &  y_2   \\
\mathrm{s.t.} & y_1 \bpx 1 & 0 & 0 \\ 0 & 0 & 0 \\ 0 & 0 & 0 \epx + y_2 \bpx 0 & 0 & 1 \\ 0 & 1 & 0 \\ 1 & 0 & 0 \epx \preceq 
\bpx 1 & 0 & 0 \\ 0 & 1 & 0 \\ 0 & 0 & 0 \epx.
\ena\right. 
\eeq
We claim that the duality gap between them is $1.$ Indeed, let $X$ be a feasible solution of (\ref{posgap}). Since $x_{11}=0, \,$ the first row and column of $X$ must be zero,  hence 
$$
X = \bpx 0 & 0 & 0 \\ 0 & 1 & 0 \\ 0 & 0 & 0 \epx  
$$
is an optimal solution with objective value $1.$ 
In turn, in (\ref{posgap-dual}) we have $y_2 = 0$ for all feasible $y,$ so its optimal value is $0.$  
  
Next,  let $\epsilon > 0$ be small, and define $M_\epsilon > 0$ 
so that 
$$
X_\epsilon := \bpx \epsilon & 0 & (1-\epsilon)/2  \\ 0 & \epsilon & 0 \\ (1-\epsilon)/2 & 0 & M_\epsilon \epx  % \, {\rm is \, positive \, semidefinite}.  
$$
is positive semidefinite. 
Then $X_\epsilon$ is an approximate solution of (\ref{posgap}), which violates only the first constraint (by $\epsilon$) 
and has  objective value $2 \epsilon.$ 

Do such ``fake" solutions arise in practice? At first look it seems that they do not.
If we  feed the pair (\ref{posgap})-(\ref{posgap-dual}) to Mosek, it returns a solution 
with DIMACS errors
$$
(0.5000, 0, 0.7071, 0, -5.5673 \cdot 10 ^{-9}, 5.9077 \cdot 10^{-17}).
$$
Since the first and third errors are large, we cannot conclude that the problem has been ``solved".

However, let us apply a similarity transformation $T^{\top}(\cdot)T$ to all matrices in (\ref{posgap}) with 
$$
T =  \bpx 3 &    5 &   -2 \\
4 &    1  &    1 \\
-4  &    -4  &    5
\epx.
$$
Then the resulting primal-dual pair still has a duality gap of $1. \,$ 
Yet, Mosek now returns a solution with DIMACS errors 
$$
(1.6093 \cdot 10^{-6},\ 0,\ 5.2111 \cdot 10^{-9},\ 3.287 \cdot 10^{-12},\ -8.1484 \cdot 10^{-5},\ 3.0511 \cdot 10^{-5}), 
$$
which may seem ``essentially all zero'' to a user. 
\eex

We argue that  in any SDP pair with positive duality gap such ``fake" solutions can arise. 
Indeed, suppose 
$$
\val\eqref{d} < \val\eqref{p}, 
$$
where $\val(\cdot)$ denotes  the optimal value of an optimization problem. 
Then by the theory of asymptotic duality (see e.g., Section 3 in \cite{renegar2001mathematical})  
there is a  sequence $\{ X_\eps \succeq 0 \mid \eps > 0 \}$ such that $X_\eps$ 
violates each primal constraint by at most $\eps, \, $ and 
$$
C \bullet X_\eps \rightarrow \val\eqref{d}  \, \, {\rm as} \,\, \eps \searrow 0.
$$
As Example \ref{ex2} shows, such ``fake" or approximate  
solutions are sometimes indeed found by SDP solvers.

%Such ``fake" solutions can arise in any SDP pair with positive duality gap. % there are such ``fake" solutions.  % indeed arise  in {\em any} 
% To explain how, let us denote by $\val(\cdot)$ the optimal value of an optimization problem, and let us  assume for simplicity 
%that  \eqref{p} and \eqref{d} are both feasible, and   
%$$
%\val\eqref{d} < \val\eqref{p}. 
%$$
%Then by the theory of asymptotic duality (see e.g., Section 3 in \cite{renegar2001mathematical})  
%there is a  sequence $\{ X_\eps \succeq 0 \mid \eps > 0 \}$ such that $X_\eps$ 
%violates each one of the primal constraints by at most $\eps, \, $ and 
%$$
%C \bullet X_\eps \rightarrow \val(D) \, \, {\rm as} \,\, \eps \searrow 0.
%$$
%As our previous example shows, such ``fake" or approximate  
 %solutions are sometimes indeed found by SDP solvers.

We note that \cite{CheWolkSchurr:12} also presented computational results on SDPs with positive duality 
gaps, and noted that Sedumi  often gave an incorrect solution on such problems. However, \cite{CheWolkSchurr:12} did not report  the DIMACS errors.

%%%%%%%
% \section{Detailed comments on some problems}
%%%%%%%

%\section{Detailed comments: how do the preprocessors perform on some of the problems?}\label{section:detailed}
\section{Detailed comments on some of the preprocessing results}\label{section:detailed}

We  now  report in detail how the  preprocessors  perform on some of the
problems. We thus examine them from several angles: for example, can they help to find known optimal solutions of difficult SDPs?  How do they perform on large-scale SDPs? 
How fast are they when they do not reduce an SDP by much, or at all? 

We first look at how the preprocessors perform on the ``Compact", ``unbound'',  and ``Example" problems, for which  the exact optimal  values  are 
known, but are hard to compute. (These problems are from the PP dataset). 
We examine whether preprocessing  helps to find  
these optimal values. % can be recovered by the preprocessing methods.

First we note that Sieve-SDP does not change the  optimal value of \eqref{p},
since it  deletes rows and columns from the variable matrix $X$ that   are always zero anyway. 
However, it deletes  %may increase the dual optimal value
rows and columns in the constraint matrices, so 
after applying it,  in the dual \eqref{d}  we require only a principal minor 
of $C - \sum_{i=1}^m y_i A_i$ to be psd.
Thus applying Sieve-SDP may increase the optimal value of (\ref{d}).

 %need to discuss what we mean  by ``recovering" the optimal values.
%the optimal value of the primal \eqref{p} does not change if we preprocess it by Sieve-SDP, but the optimal value of the dual \eqref{d} may increase.
%Indeed, the primal optimal value stays the same, since Sieve-SDP  
%only deletes rows and columns from $X$ which  are always zero anyway. 
%As to the dual optimal value,  
%since Sieve-SDP  deletes  rows and columns in the constraint matrices, 
%after preprocessing we only require a principal minor 
%of $C - \sum_{i=1}^m y_i A_i$ to be psd.
%Thus the optimal value of (\ref{d}) may indeed  increase.

%First we need to discuss what we mean  by ``recovering" the optimal values. %a known optimal solution. 
%When we preprocess by Sieve-SDP the optimal value of the primal \eqref{p} does not change, 
%since it only deletes rows and columns from $X$ which  are always zero anyway. 
%	As to the dual optimal value,  
	%	since Sieve-SDP  deletes  rows and columns in the constraint matrices, 
	%after preprocessing we only require a principal minor 
	%of $C - \sum_{i=1}^m y_i A_i$ to be psd.
	%Thus the optimal value of (\ref{d}) may   increase.
	
To quantify %and illustrate 
this argument, let 
	$(P_\text{pre})$ and $(D_\text{pre})$ be the primal and dual problems after preprocessing by Sieve-SDP, respectively.
	Then  %by the previous argument 
	\beq \label{ineq-sieve} 
	\val \eqref{d}  \leq \val(D_\text{pre}) \leq \val(P_\text{pre}) =  \val\eqref{p}. 
	\eeq
First assume $	\val \eqref{d}  < \val\eqref{p}.$ Then we can show by examples that  any  inequality  in 
	\eqref{ineq-sieve} 	may be strict. For example, in Example \ref{ex2}   Sieve-SDP deletes the first row and first column in all constraint matrices,
	and it is easy to check that the corresponding optimal values are $0 < 1 = 1 = 1, \,$ respectively. In detail, for this example 
	   $(D_\text{pre})$   is 
	\beq \label{posgap-dual-pre-1} 
	\left.\ba{rl}
	\displaystyle\sup_{y_2}   &  y_2   \\
	\mathrm{s.t.} &  y_2 \bpx  1 & 0 \\  0 & 0 \epx \preceq   \bpx   1 & 0 \\  0 & 0 \epx,
	\ena\right.
	\eeq
	whose optimal value is $1.$
	
On the other hand,  suppose  $\val \eqref{p} = \val \eqref{d}. \,$ Then in \eqref{ineq-sieve}  equality holds throughout, so  Sieve-SDP changes neither the  primal, nor the dual optimal values. 

%Of course, \eqref{ineq-sieve} implies that when $\val \eqref{p} = \val \eqref{d}, \, $ then preprocessing by Sieve-SDP changes neither t he primal, nor the dual optimal values. 
	
Which optimal values are changed or kept the same by the other preprocessors? 	Pd1 and pd2 also reduce the primal \eqref{p}, so when we apply them, % by them, 
the primal optimal value (but maybe not that of the dual) will remain the same.
	On the other hand,  dd1 and dd2  reduce the dual problem (\ref{d}), 
	so they keep its optimal value  the same. However,   they may change the optimal value of 
	the primal \eqref{p}.

In all tables in this section 
we use the following convention: the first reported objective value is the primal and the second is the dual.

% \subsection{``Compact" problems -- 10 problems from  \cite{Waki:12}}

\subsection{``Compact" problems -- 10 problems from  \cite{Waki:12}}\label{subsect:compact}

These instances  are   {\em weakly infeasible}, i.e., 
the affine subspace 
$$
H \, = \, \{ \, X \, \mid \, A_i \bullet X = b_i \, (i=1, \cdots, m) \}
$$ 
does not intersect $\psd{n},$ but the distance of $H$ to $\psd{n}$ is zero. Weakly infeasible SDPs are particularly challenging to SDP solvers. However, 
a recent algorithm in \cite{HenrionNaldi:2016} can detect (in)feasibility of small SDPs in exact arithmetic,  and 
  \cite{Liu2018}  presented an algorithm that is tailored to detect weak infeasibility.

On these problems pd1 and pd2 produced the same results, 
%Pd1 and pd2 produced the same results, 
while dd1 and dd2 reduced none of them. Pd1 and pd2 combined with Mosek correctly detected primal infeasibility of all problems, 
while Sieve-SDP correctly proved  primal   infeasibility without  Mosek.
(Since it found the primal infeasible, we did not compute a dual solution).

The results are in Table \ref{table:compactdim}.

\begin{table}[!htbp]
\newcommand{\crf}[1]{{}#1{}}
\newcommand{\crfs}[1]{{}#1{}}
\begin{center}
\caption{Results on the ``Compact" problems}\label{table:compactdim}
\setlength{\tabcolsep}{2pt}
\resizebox{\textwidth}{!}{	
	\begin{tabular}{lccccc}
		\hline\noalign{\smallskip}
		\crf{Problem} & \crf{Correct obj (P, D)} & \crf{Obj before} & \crf{After pd1/pd2} & \crf{After dd1/dd2} & \crf{After Sieve-SDP} \\
		\noalign{\smallskip}\hline\noalign{\smallskip}
		\crf{CompactDim2R1	} & \crfs{Infeas, $+\infty$}	& 3.79e+06, 4.20e+06	& Infeas, 1	& 3.79e+06, 4.20e+06& Infeas, - \\
		\crf{CompactDim2R2	} & \crfs{Infeas, $+\infty$}	& 6.41e-10, 6.81e-10	& Infeas, 2	& 6.41e-10, 6.81e-10 & Infeas, - \\
		\crf{CompactDim2R3	} & \crfs{Infeas, $+\infty$}	& 1.5, 1.5				& Infeas, 2	& 1.5, 1.5			& Infeas, - \\
		\crf{CompactDim2R4} & \crfs{Infeas, $+\infty$}	& 1.5, 1.5				& Infeas, 2	& 1.5, 1.5			& Infeas, - \\
		\crf{CompactDim2R5	} & \crfs{Infeas, $+\infty$}	& 1.5, 1.5				& Infeas, 2	& 1.5, 1.5			& Infeas, - \\
		\crf{CompactDim2R6} & \crfs{Infeas, $+\infty$}	& 1.5, 1.5				& Infeas, 2	& 1.5, 1.5			& Infeas, - \\
		\crf{CompactDim2R7	} & \crfs{Infeas, $+\infty$}	& 1.5, 1.5				& Infeas, 2	& 1.5, 1.5			& Infeas, - \\
		\crf{CompactDim2R8	} & \crfs{Infeas, $+\infty$}	& 1.5, 1.5				& Infeas, 2	& 1.5, 1.5			& Infeas, - \\
		\crf{CompactDim2R9	} & \crfs{Infeas, $+\infty$}	& 1.5, 1.5				& Infeas, 2	& 1.5, 1.5			& Infeas, - \\
		\crf{CompactDim2R10} & \crfs{Infeas, $+\infty$}	& 1.5, 1.5				& Infeas, 2	& 1.5, 1.5			& Infeas, - \\
		\noalign{\smallskip}\hline\noalign{\smallskip}
		Correctness  \% & 100\%, 100\% & 0\%, 0\% & 100\%, 0\% & 0\%, 0\% & 100\%, - \\
		\noalign{\smallskip}\hline
	\end{tabular}
} 
\end{center}
\vspace{-3ex}
\end{table}

We mention here  another set of  infeasible, and weakly infeasible 
SDPs. They are described  in 
	\cite{LiuPataki:17}, and are available from the webpage of G\'abor Pataki. 
Some of these SDPs are classified as ``clean" and some of them  as 
``messy".  In the ``clean" instances the structure that proves infeasibility is apparent, while 
in the ``messy" instances that structure was obscured by two kinds of operations: random elementary row operations on the constraints and a 
random similarity transformation.

Indeed, in our testing  all clean instances were found infeasible 
by  Sieve-SDP, pd1,  and pd2.  %  detected infeasibility of
 In contrast, no messy instances were reduced by any of the preprocessors. 
Since the clean instances are evidently easy for Sieve-SDP, and the messy ones are hard for all preprocessors, 
we did not include the SDPs from 	\cite{LiuPataki:17} in our test set, since we felt that this would not be fair. 

% \subsection{``unbound'' problems -- 10 problems from \cite{waki2012strange}}

\subsection{``unbound'' problems -- 10 problems from \cite{waki2012strange}}\label{subsect:unbounded}

The mathematically correct optimal values of both the primal and the dual are $0$  in this problem collection.
%  the  optimal values are $0$ for both the primal and dual problems.
However, before preprocessing Mosek returned wrong optimal  values for $6$ out of $10$ problems.  
Although Mosek found  solutions with almost correct optimal  value in  problems 2, 3 and 4, these solutions are inaccurate, as 
the DIMACS errors are of the order $10^{-1}$ (this is marked by "*" symbols in Table \ref{table:unbounddim}).

In summary, $9$ out of $10$ 
problems in this dataset need preprocessing to obtain  a reasonable solution.

Sieve-SDP, pd1 and pd2 corrected  all objective values,  as Table  \ref{table:unbounddim} shows.

It is interesting that the authors in \cite{waki2012strange} computed the correct optimal solution of  these instances using SDPA-GMP \cite{fujisawa2008sdpaGMP}, a high-precision SDP solver that carries several hundred significant digits. Of course, running SDPA-GMP  is more time consuming, than running Sieve-SDP and Mosek.

\begin{table}[!htbp]
%\vspace{-2ex}
\begin{center}
	\caption{Results on the ``unbound" problems}\label{table:unbounddim}
	\setlength{\tabcolsep}{2pt}
	\resizebox{\textwidth}{!}{		
	\begin{tabular}{lccccc}
		\hline\noalign{\smallskip}
		Problem & Correct obj (P, D) & Obj before & After pd1/pd2 & After dd1/dd2 & After Sieve-SDP \\
		\noalign{\smallskip}\hline\noalign{\smallskip}
		unboundDim1R1	& 0, 0	& 1.33e-09, -7.05e-10 		& 1.33e-09, -7.05e-10 & 1.33e-09, -7.05e-10 & 0, 0 \\
		unboundDim1R2	& 0, 0	& -8.19e-15*, -8.01e-15*	& 0, 0	& -8.19e-15*, -8.01e-15*	& 0, 0 \\
		unboundDim1R3	& 0, 0	& -2.04e-11*, -2.02e-11*	& 0, 0	& -2.04e-11*, -2.02e-11*	& 0, 0 \\
		unboundDim1R4	& 0, 0	& -2.34e-10*, -2.32e-10*	& 0, 0	& -2.34e-10*, -2.32e-10*	& 0, 0 \\
		unboundDim1R5	& 0, 0	& -1, -1	& 0, 0	& -1, -1	& 0, 0 \\
		unboundDim1R6	& 0, 0	& -1, -1	& 0, 0	& -1, -1	& 0, 0 \\
		unboundDim1R7	& 0, 0	& -1, -1	& 0, 0	& -1, -1	& 0, 0 \\
		unboundDim1R8	& 0, 0	& -1, -1	& 0, 0	& -1, -1	& 0, 0 \\
		unboundDim1R9	& 0, 0	& -1, -1	& 0, 0	& -1, -1	& 0, 0 \\
		unboundDim1R10	& 0, 0	& -1, -1	& 0, 0	& -1, -1	& 0, 0 \\
		\noalign{\smallskip}\hline\noalign{\smallskip}
		Correctness \% & 100\%, 100\% & 10\%, 10\% & 100\%, 100\% & 10\%, 10\% & 100\%, 100\% \\
		\noalign{\smallskip}\hline
	\end{tabular}
	}
\end{center}	
\vspace{-3ex}	
\end{table}

% \subsection{``Example" problems -- 8 problems from \cite{CheWolkSchurr:12}}

\subsection{``Example" problems -- 8 problems from \cite{CheWolkSchurr:12}}\label{subsect:example}

The mathematically correct objective values are reported in \cite{CheWolkSchurr:12} in table 12.1.
(Note that in \cite{CheWolkSchurr:12} our primal is considered  the  dual, and vice versa, so that table must be read accordingly.)

Table \ref{table:example} shows the objective values before and after preprocessing. We consider an objective value correct if it is less than ${10}^{-6}$ away from the true optimal value.

We excluded ``Example5'' of \cite{CheWolkSchurr:12} 
 from this table, since in Table 12.1 in \cite{CheWolkSchurr:12} its optimal value is not reported. For all other problems, except for ``Example9size20'' and ``Example9size100'', we manually verified the correctness of the optimal values in exact arithmetic.

\begin{table}[!htbp] 
%\vspace{-2ex}
\begin{center}
	\caption{Results on the ``Example" problems}\label{table:example}
	\setlength{\tabcolsep}{2pt}
	\resizebox{\textwidth}{!}{		
	\begin{tabular}{lccccc}
		\hline\noalign{\smallskip}
		Problem & Correct obj (P, D) & Obj before & After pd1/pd2 & After dd1/dd2 & After Sieve-SDP \\
		\noalign{\smallskip}\hline\noalign{\smallskip}
		Example1		& 0, 0			& 0, 0				& 0, 0				& 0, 0				& 0, 0 \\
		Example2		& 1, 0			& 3.33e-01, 3.33e-01& 1, 1				& 4.73e-15, 1.82e-14& 1, 1 \\
		Example3		& 0, 0			& 3.33e-01, 3.33e-01& 1.17e-07, 1.69e-07& 4.73e-15, 1.82e-14& 1.17e-07, 1.69e-07 \\
		Example4		& Infeas, 0{~~~~~~}		& {~~~}Infeas, 3.74e-07	& Infeas, 1{~~~~~~}			& 0, 0				& Infeas, -{~~~~~~}	 \\
		Example6		& 1, 1			& 1, 1				& 1, 1				& 1, 1				& 1, 1 \\
		Example7		& 0, 0			& 0, 0				& 0, 0				& 0, 0				& 0, 0 \\
		Example9size20 	& Infeas, 0{~~~~~~}		& {~~~}Infeas, 3.39e-01	& Infeas, 1{~~~~~~}			& 0, 0				& Infeas, -{~~~~~~}	\\
		Example9size100 & Infeas, 0{~~~~~~}		& {~~~}Infeas, 3.43e-01	& Infeas, 1{~~~~~~}			& 0, 0				& Infeas, -{~~~~~~} \\
		\noalign{\smallskip}\hline\noalign{\smallskip}
		Correctness  \%	& 100\%, 100\%	& 75\%, 50\%		& 100\%, 50\%		& 50\%, 100\%		& 100\%, 50\% \\
		\noalign{\smallskip}\hline
	\end{tabular}}
\end{center}
\vspace{-3ex}	
\end{table}

Note that the comparison in Table \ref{table:example} is somewhat unfair to Sieve-SDP:
if it found a problem infeasible, it did not compute a dual solution.

% \subsection{``finance'' problems -- 4 problems from \cite{boyd2014performance}}

\subsection{``finance'' problems -- 4 problems from \cite{boyd2014performance}}\label{subsection-finance}

The PP dataset contains four ``finance" problems: ``leverage\_limit'', ``long\_only'', ``sector\_neutral'' and ``unconstrained''.
We report on these problems in detail, since these are the largest 
in the PP dataset. 
For example, ``long\_only'' has $100$ semidefinite variable blocks of order $91$ and another 
$100$ of order $30.$

%In the PP dataset these problems are the largest, so we report on them in detail.
%For example, the problem ``long\_only'' has $100$ semidefinite variable blocks of order $91$ and another 
%$100$ of order $30.$

Table \ref{table:finance} shows how much the preprocessors reduced these SDPs:  here $n_{\textrm{sdp}}$  is the total size of the semidefinite  blocks; 
$n_{\textrm{nonneg}}$ is the total number of nonnegative variables; 
$n_{\textrm{free}}$  is the total number of free variables;  $m$ is the total number of constraints; and $\mathrm{nnz}$ is the total number of nonzeros.

While  dd1 and dd2 significantly reduced the size of the SDP blocks, they added many free variables. 
Sieve-SDP reduced the size of the SDP blocks, without adding free variables, and it eliminated the most constraints. 
We mention that after preprocessing with dd2 Mosek detected  that  problem 
 ``leverage\_limit'' is ``dual infeasible.''  
This may be because of numerical instability, and does not contradict the result we get after preprocessing with Sieve-SDP.

\begin{table}[!htbp]
%\vspace{-2ex}
\begin{center}
	\caption{Results on the ``finance'' problems}\label{table:finance}
	%\resizebox{\textwidth}{!}{	
	\begin{tabular}{lrrrrrrr}
		\hline\noalign{\smallskip} Method 
		& $n_{\textrm{sdp}}$ & $n_{\textrm{nonneg}}$ & $n_{\textrm{free}}$ & $m$ & $\mathrm{nnz}$\\
		\noalign{\smallskip}\hline\noalign{\smallskip}
		None 		& 60,400	& 51,100	& 0	& 251,777
& 2,895,756\\
		After pd1	& 60,400	& 51,100	& 0	& 251,777
& 2,895,756\\
		After pd2	& 60,280	& 51,100	& 0	& 249,797
& 2,880,876\\
		After dd1	& 27,429	& 51,100	& 2,286,000	& 251,777
& 2,844,756\\
		After dd2	& 36,400	& 51,100	& 2,521,005	& 251,777
& 2,605,807\\
		After Sieve-SDP & 56,766 & 50,873	& 0 & 215,210
& 2,466,573\\
		\noalign{\smallskip}\hline
	\end{tabular}%}
\end{center}
\vspace{-3ex}	
\end{table}

We remark that preprocessing actually {\em increased} the solution time on these problems, though not by much.  
For example, 
the total time spent on preprocessing with Sieve-SDP plus solving with Mosek is about 21\% higher than the solving time with Mosek without preprocessing. Still, since the primary goal of preprocessing is 
to improve solution accuracy, we believe that we should do it whenever we can.

Furthermore, on these instances Sieve-SDP performed a large number of iterations, and deleted only a small submatrix in each one. Thus, we could easily reduce the time spent by Sieve-SDP by limiting the maximum number of iterations it is allowed to perform. 
We do not report results with such a setting, since we do not want to ``overtune" our code.

% \subsection{Dressler-Illiman-de Wolff (DIW) dataset (155 problems)}

\subsection{Dressler-Illiman-de Wolff (DIW)  dataset (155 problems)}

Consider the optimization problem
\begin{equation} \label{problem-polyopt} 
\left.\begin{array}{ll}
\displaystyle\min_x & f(x)\\
\mathrm{s.t.}	& g_i(x) \geq 0~~ (i=1, \dots, m), 
\end{array}\right.
\end{equation} 
where $f$ and the $g_i$ are multivariate polynomials.
As shown in the seminal work of Lasserre \cite{lasserre2001global}, 
the optimal value of (\ref{problem-polyopt}) 
can be lower bounded by solving  SDPs. 
Under suitable conditions the  lower bounds converge to the 
optimal value of (\ref{problem-polyopt}),  as the so-called Lasserre relaxation 
order increases.
However, no useful lower bound is obtained when the SDPs are infeasible. See Parrilo \cite{parrilo2003semidefinite} for 
a related scheme to construct  SDP relaxations of  (\ref{problem-polyopt}). 

Since solving the Lasserre SDPs can be challenging, Dressler, Illiman and de Wolff  \cite{dressler2018approach} 
proposed an alternative relaxation, based on so-called nonnegative circuit polynomials, 
and they compared their approach with the SDP-based one.

We constructed  the SDPs in the ``DIW" dataset by taking the polynomial optimization problems from 
\cite{dressler2018approach} and using  Gloptipoly 3 (\cite{henrion2009gloptipoly}) to generate their SDP relaxations.

We describe our SDPs in Table~\ref{tbl:dressler_info} 
 with their Lasserre relaxation order, which % of the SDPs is also shown.
 ranges  
from the lowest possible (half the degree of the highest degree monomial in the polynomials) to 20.  
For example, the SDP named ``ex3.3\_order4''  is obtained  by applying the Lasserre relaxation of order 4 to Example 3.3 in 
\cite{dressler2018approach}.

\begin{table}[!htbp]
%\vspace{-2ex}
\begin{center}
	\caption{Relaxation orders for examples in \cite{dressler2018approach} \label{tbl:dressler_info}}
	\setlength{\tabcolsep}{2pt}
	\resizebox{\textwidth}{!}{	
	\begin{tabular}{l c c c c c c c c c }
		\hline\noalign{\smallskip}
		ex 	& $3.3$		& $4.1$		& $4.2$		& $4.3$		& $4.4$		& $5.4$		& $5.5$		& $5.6$		& $5.7$\\
		\noalign{\smallskip}\hline\noalign{\smallskip}
			relaxation orders	& $6 \cdots 20$	 & $3 \cdots 20$	& $6 \cdots 20$	& $2 \cdots 20$	& $3 \cdots 20$	 & $5 \cdots 20$	& $4 \cdots 20$	 & $4 \cdots 20$	& $5 \cdots 20$ \\
		\noalign{\smallskip}\hline
	\end{tabular}}
\end{center}	
\vspace{-3ex}
\end{table}

Table \ref{tbl:result_of_dressler} shows the  results:  ``$n$'' is the sum of the orders of all psd and nonnegative blocks, and 
	``$m$" is the sum of the number of constraints in all problems.
	
The results are quite striking.
\begin{table}[!htbp]
	%\vspace{-2ex}
	\begin{center}
		\caption{Results for the DIW dataset}\label{tbl:result_of_dressler}
		%\setlength{\tabcolsep}{2pt}
		%\resizebox{\textwidth}{!}{	
			\renewcommand{\arraystretch}{1.2}
			\begin{tabular}{l | r  r r  | r r | r } \hline
				%		\hline\noalign{\smallskip}
				\rule{0pt}{12pt}Method & \# Reduced & $n$ & $m$ & Preprocessing (s) & Solving (s) & \# Infeas\\ \hline
				%		\noalign{\smallskip}\hline\noalign{\smallskip}
				None		& -  {~}	& 53,523 & 186,225 & -  {~}  & 139,493.56 & -  {~} \\
				pd1			& 155 {~}	& 1,450	& 3,278		& 1632.43  {~}	& 128.46	& 56\\
				pd2			& 155 {~}	& 1,450	& 3,278		& 10,831.32 {~}	& 124.44	& 56\\
				dd1			& 0 {~}		& 53,523	& 186,225	& 65.18 {~}		& 139,493.56	& 0\\
				dd2			& 0 {~}		& 53,523	& 186,225	& 22,152.57 {~}	& 139,493.56	& 0\\
				Sieve-SDP	& 155 {~}	& 1,385 & 3,204 & 1,232.27 {~} & 87.53 & 59\\ \hline
				%		\noalign{\smallskip}\hline
			\end{tabular}%}
		\end{center}
		\vspace{-3ex}
	\end{table}
Sieve-SDP, pd1, and pd2  %very well: they
 ran fast, reduced all problems, detected infeasibility of more than a third, and 
reduced overall computing time by a factor of more than a hundred! 
Sieve-SDP was the best in all aspects, with pd1 a close second.
%  with respect to both preprocessing time and solving time after preprocessing. Pd1 was a close second. 

Note that without preprocessing Mosek failed  to detect infeasibility of any of these SDPs.

These results are somewhat surprising since \cite{dressler2018approach} solved some of these SDPs to near optimality, and  managed to extract approximate  optimal solutions of the original polynomial optimization problems.
See  \cite{henrion2005detecting} for similar results on  similar SDPs. % were obtained  earlier in \cite{henrion2005detecting}. 
In fact,  \cite{henrion2005detecting}  took the view that numerical inaccuracy of the SDP solvers  actually {\em helps} find 
near-optimal solutions 
of the polynomial optimization problems. See \cite{lasserre2018sdp} for a more recent and thorough study of the same issue.

We remark that these SDPs are likely to be weakly infeasible.

We were thus motivated to  
double check  that Sieve-SDP indeed reduced these SDPs 
  	correctly.  Precisely, we verified that in the Basic Step (in Figure \ref{figure-basicstep}) it  only eliminated  constraints 
  in one of the following forms: either of the form
  	\beq \nonumber
  	\bpx D & 0 \\ 0 & 0 \epx \bullet X = 0, 
  	\eeq
  	where $D$ is positive definite diagonal, of order $1$ or $2, \,$ and 
  		the smallest diagonal element 
  		is  $1$ or $0.5$ or $1/3=0.3333...;$ or of the form 
  		$$O \bullet  X = 0,$$
  		where $O$ is the zero matrix. 
  		Furthermore,  Sieve-SDP  always detected infeasibility  by finding a constraint 
  		\beq \nonumber 
  		\bpx D  & 0 \\ 0 & 0 \epx \bullet X = \beta,  
  		\eeq
  		where $D$ is as above, and $\beta = -3$ or $-8.$

The zeroes  in all these constraints are  zeroes  in absolute machine precision, i.e., in the sparse SDPs returned by Gloptipoly 3
  	these entries do not appear at all. 
  	Thus Sieve-SDP performed all reductions correctly.

% \subsection{Henrion-Toh  dataset (98 problems)}

\subsection{Henrion-Toh  dataset (98 problems)}

This dataset  was 
kindly provided to us by Didier Henrion and Kim-Chuan Toh.
The problems come mostly from polynomial optimization.

Among these problems $18$ were reduced by pd1, pd2, or Sieve-SDP and none 
	by dd1 or dd2. 
	Table  \ref{table-henrion}  shows the time details in seconds.
	The last column ``Pre. vs. Solve" shows the time spent on preprocessing as a percentage of time spent on solving. It is
	\begin{equation} \label{eqn-pre-vs-solve} 
	\dfrac{\rm preprocessing \, time}{\rm solving  \, time \, without \, preprocessing} \times 100 \, \%. 
	\end{equation}
	  
	\begin{table}[!htbp]  
		%\vspace{-2ex}
		\begin{center}
			\caption{Time results on the  Henrion-Toh dataset}\label{table-henrion}
			%\resizebox{\textwidth}{!}{	
				\begin{tabular}{l r r r r}
					\hline\noalign{\smallskip}
					Method & Preprocessing (s) & Solving (s) & Pre.  vs. solve  \\
					\noalign{\smallskip}\hline\noalign{\smallskip}
					None		& - 	& 1420.02	& -			\\
					pd1			& 10.27	& 1373.70	& 0.72\%	\\
					pd2 		& 49.84 & 1374.31	& 3.51\%	\\
					dd1			& 3.93	& 1420.02	& 0.28\%	\\
					dd2			& 29.24	& 1420.02	& 2.06\% \\
					Sieve-SDP	& 4.58	& 1376.27	& 0.32\% \\
					\noalign{\smallskip}\hline
				\end{tabular}%}
		\end{center}
		\vspace{-3ex}	
	\end{table}

On this dataset  the  preprocessors are less successful: 
pd1, pd2, and Sieve-SDP detected infeasibility of only one problem (of ``sedumi-l4")
and they reduced solving time only a little. % moderately. %by much.
%they did not detect infeasibility in any one of the instances,  and 
%the time reduction is not significant. 
However,  the  preprocessing times are small, or even negligible: 
for example, Sieve-SDP spent only about $0.3 \%$ of the time that it took for Mosek to solve the problems.

\begin{figure}[!htbp]%[H]
	\begin{center}
		\includegraphics[width = 7cm, height = 7cm]{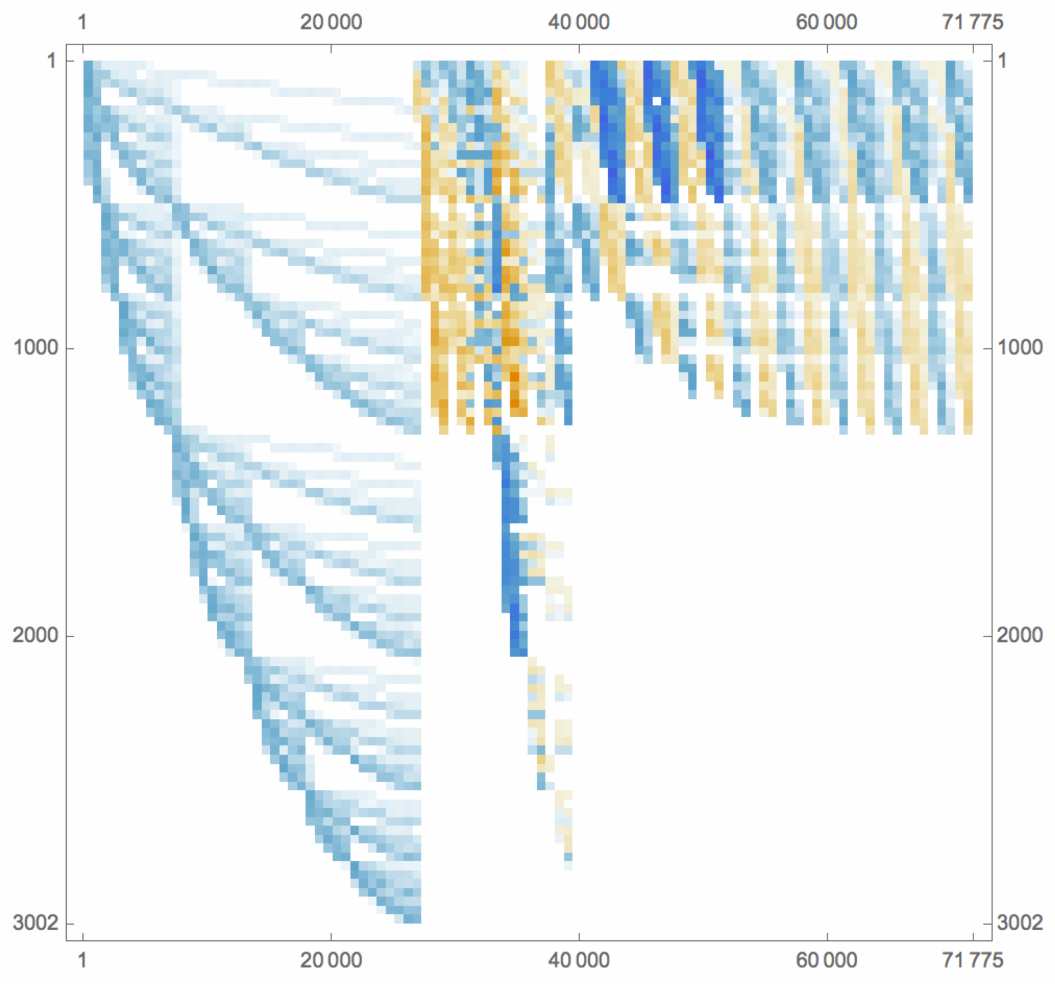} \hspace{0.1cm}\includegraphics[width = 4cm, height = 4cm]{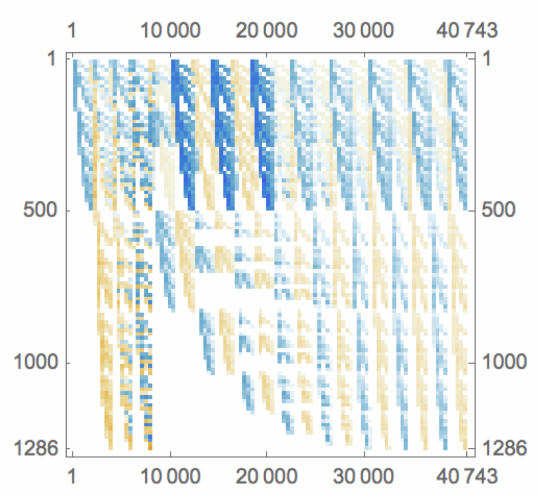} 
		\caption{Instance ``sedumi-fp32": size and sparsity before (left) and after (right) preprocessing}\label{figure:sedumifp32} 
	\end{center}	
\vspace{-3ex}	
\end{figure}

In Figure \ref{figure:sedumifp32} we illustrate how Sieve-SDP works on the instance ``sedumi-fp32'': we show the sparsity structure of the constraints of  the original problem (on the left),
and after Sieve-SDP (on the right).   Just like in Figure \ref{figure:DIW}, 
each row corresponds to an $A_i$ matrix stretched out as a vector. Red dots correspond to positive entries, blue dots  correspond to negative entries, and white areas to zero entries.

Here we also discuss problem ``sedumi-fp33'' on which preprocessing by  
Sieve-SDP makes the DIMACS error {\em worse}.
Since this is the only such instance, we looked at it in more detail. The worst DIMACS error (of a solution computed by Mosek) before Sieve-SDP is $3.36\times{10}^{-7},$ which is acceptable. After Sieve-SDP the worst error is about $0.0928, \,$ which is unacceptable.

We also solved this instance using the high accuracy SDP solver SDPA-GMP \cite{fujisawa2008sdpaGMP}. 
The DIMACS errors were 
$$
2.3497 \cdot 10^2, 0.0000, 1.8552 \cdot 10^1, 0.0000, -9.9999 \cdot 10^{-1}, 8.5173 \cdot 10^{-2} 
$$
before Sieve-SDP, and 
$$
3.4075 \cdot 10^2, 0.0000, 1.9636 \cdot 10^1, 0.0000, -9.9999 \cdot 10^{-1}, 6.1901 \cdot 10^{-1}
$$
after Sieve-SDP. In both cases the largest error is more than $200, \,$ which is unacceptably large.

Given the high accuracy of SDPA-GMP,  it seems  that 
this SDP cannot be accurately solved by current fast solvers, and 
the {\em worse}  DIMACS error returned by Mosek after Sieve-SDP alerts the user to this fact:
this problem may actually have a positive duality gap 
(cf.  Example \ref{ex2}).

% \subsection{Toh-Sun-Yang dataset (419 problems)}

\subsection{Toh-Sun-Yang dataset (419 problems) from \cite{sun2015convergent, yang2015sdpnal}}

Although none of the five methods reduced the SDPs in this collection, we still comment on them in detail, since 
 it is interesting that 
pd1, dd1 and Sieve-SDP spent only a negligible amount time on preprocessing. Thus using these 
 three methods 
it does not hurt to preprocess: see 
Table \ref{tbl:TSY_collection}. The last column ``Pre. vs. Solve" shows the time spent on preprocessing as a percentage of time spent on solving; see 
equation \eqref{eqn-pre-vs-solve}. Pd2 and dd2, on the other hand, spent considerably more time on preprocessing.
 
 \begin{table}[!htbp]
 	%\vspace{-2ex}
 	\begin{center}
 		\caption{Timing on the Toh-Sun-Yang dataset}\label{tbl:TSY_collection}
 		%\begin{tabular}{l | r r r | r }
 		\begin{tabular}{l r r r }
 			\hline\noalign{\smallskip}
 			Method & Preprocessing & Solving & Pre. vs. solve       \\
 			\noalign{\smallskip}\hline\noalign{\smallskip}
 			pd1			& 220.18	& 27,635.46	& 0.80\%\\
 			pd2			& 4,029.61	& 27,635.46	& 14.58\%\\
 			dd1			& 134.64	& 27,635.46	& 0.49\%\\
 			dd2			& 2,428.82	& 27,635.46	& 8.79\%\\
 			Sieve-SDP	& 152.14	& 27,635.46	& 0.55\%\\
 			\noalign{\smallskip}\hline
 		\end{tabular}
 	\end{center}	
 	\vspace{-3ex}
 \end{table}
 
%%%%%%%
% \section{Summary}
%%%%%%%

\section{Summary}\label{section:summary}

We now compare all preprocessors on all instances in Tables \ref{table-redandhelp}, 
\ref{table:times},  and \ref{table-redsize}.

In Table \ref{table-redandhelp} the second  column shows how many problems were reduced.
The third column shows how many problems were found to be infeasible.
The fourth column shows on how many instances the preprocessing improved the 
DIMACS errors, as we discussed in Subsection \ref{subsection-criteria-comparison}.

The last column ``Memory" shows how many times a method   ran out of memory, or crashed: 
	this happened with  pd2 six times and with dd2 four times. To ensure fair reporting 
	we reran these methods on the same instances   on a machine with 24 GB RAM, and the results were the same.

\begin{table}[!htbp]
%\vspace{-2ex}
\begin{center}
	\caption{Infeasibility detection and error reduction on all 771 problems}\label{table-redandhelp}
	\begin{tabular}{l r r r r}
		\hline\noalign{\smallskip}
		\multirow{2}{*}{Method} & \multirow{2}{*}{ \# Reduced} & \# Infeas  & \# DIMACS  error & \multirow{2}{*}{ Memory} \\
		  & & detected & improved &        \\
		\noalign{\smallskip}\hline\noalign{\smallskip}
		pd1			& 209	& 67	& 74	& 0\\
		pd2			& 230	& 67	& 78	& 6\\
		dd1			& 14	& 0		& 2		& 0\\
		dd2			& 21	& 0	 	& 4  	& 4\\
		Sieve-SDP	& 216	& 73 	& 74	& 0\\
		\noalign{\smallskip}\hline
	\end{tabular}
\end{center}	
\vspace{-3ex}
\end{table}

Table \ref{table:times} shows the preprocessing and solving times in seconds. 
The first column shows the preprocessing time 
and the second shows the solving time by Mosek after preprocessing.
Column ``Pre vs. solve" shows the relative speed of the preprocessors: see equation \eqref{eqn-pre-vs-solve}.
The last column, ``Time reduction'',  shows  by how much preprocessing decreased 
	the solving time. It is  
\begin{equation*}
%{\scriptsize
	\frac{\text{solving time w.o.   preprocessing} - \text{(preprocessing time + solving time after preprocessing)}}{\text{solving time w.o.  preprocessing}} \times 100 \, \%.
%}	
\end{equation*}
	Of course, the higher this percentage, 
	the more a preprocessor reduces solution  time.
	A negative percentage  means that preprocessing  actually {\em increased}  the total time.

\begin{table}[!htbp]
%\vspace{-2ex}
\begin{center}
	\caption{Preprocessing and solving  times on all 771 problems}\label{table:times}
	\begin{tabular}{l r r r r}
		\hline\noalign{\smallskip}
		Method & Preprocessing (s) 	& Solving	(s) 	& Pre vs. solve		& Time reduction \\
		\noalign{\smallskip}\hline\noalign{\smallskip}
		none 	& -				& 272,427.23	& -			& -\\
		pd1		& 2,486.51		& 132,356.63	& 0.91\%	& 50.50\%\\
		pd2		& 23,323.07		& 131,636.47	& 8.56\%	& 43.12\%\\
		dd1		& 587.93		& 272,244.62	& 0.22\%	& -0.15\%\\
		dd2		& 35,984.45		& 272,031.04	& 13.21\%	& -13.16\%\\
		Sieve-SDP & 2,170.13	& 131,837.25	& 0.80\%	& 51.81\%\\
		\noalign{\smallskip}\hline
	\end{tabular}
\end{center}
\vspace{-3ex}	
\end{table}

Finaly, Table  \ref{table-redsize} 
shows by ``how much'' the problems were reduced. As in Table \ref{table-redandhelp}, the second column shows the number of problems reduced by each method.

To explain the other columns, let us fix an SDP in the primal form \eqref{p} with 
	potentially several semidefinite block variables (some of which may be of order $1, \,$ i.e., they may be just nonnegative variables).
	
Let $n_{\rm before} \, {\rm and} \, n_{\rm after}$ be the total size of the semidefinite blocks before and after reduction. 
	We define the reduction rate on $n$ as 
	\begin{equation*}
	\dfrac{\sum {n_{\rm before}} - \sum {n_{\rm after}}}{\sum {n_{\rm before}}}, 
	\end{equation*}
	where the sum is  over all $771$ problems.
	
Similarly, let $m_{\rm before} \, {\rm and} \, m_{\rm after}$ be the number of constraints in a problem before and after reduction. 
	We define the reduction rate on $m$ as 
	\begin{equation*}
	\dfrac{\sum {m_{\rm before}} - \sum {m_{\rm after}}}{\sum {m_{\rm before}}}, 
	\end{equation*}
	where the sum is again taken over all $771$ problems.

Methods dd1 and dd2 added  free variables, and the fifth column in Table~\ref{table-redsize}
	shows how many.

The sixth column ``nnz" shows the total number of nonzeros in the constraint matrices.

\begin{table}[!htbp]
%\vspace{-2ex}
\begin{center}
	\caption{Size reduction on all 771 problems}\label{table-redsize}
	\begin{tabular}{l r r r r r}
		\hline\noalign{\smallskip}
		Method & \# Reduced & Red. on $n$ & Red. on $m$ & Extra free vars & nnz\\
		\noalign{\smallskip}\hline\noalign{\smallskip}
		none		& -		& -			& -			& -			& 300,989,332\\
		pd1			& 209	& 15.47\%	& 17.79\%	& 0			& 211,299,702\\
		pd2			& 230	& 15.59\%	& 18.23\%	& 0			& 211,257,726\\
		dd1			& 14	& 6.74\%	& 0.00\%	& 2,293,495	& 300,936,120\\
		dd2			& 21	& 9.28\%	& 0.00\%	& 2,315,849	& 299,272,012\\
		Sieve-SDP	& 216	& 16.55\%	& 20.66\%	& 0			& 206,061,059\\
		\noalign{\smallskip}\hline
	\end{tabular}
\end{center}
\vspace{-3ex}	
\end{table}

Given these tables we  now  summarize our  findings. 
In all aspects Sieve-SDP is competitive with the other preprocessing methods. 
In detail:

\begin{itemize}
	\item It is competitive  considering the number of problems reduced. 
	\item It is competitive in computing known optimal solutions; see Tables   \ref{table:compactdim}, \ref{table:unbounddim},  and \ref{table:example}.
	\item The time spent on preprocessing with Sieve-SDP vs. solving is negligible. It is also negligible for  pd1 and  dd1, but less so for pd2 and dd2. See Table \ref{table:times}.
\end{itemize}

\noindent In several aspects Sieve-SDP is the best.

\begin{itemize}
			\item It is best in detecting infeasibility: see Table \ref{table-redandhelp}. 
		It is important that  Sieve-SDP detects infeasibility without using 
		any optimization solver, whereas the other methods 
		rely on Mosek. % (both to solve the LP subproblems, and to detect infeasibility).
		\item It  reduced solution time the most,   with pd1 a close second. See Table \ref{table:times}. 
		\item It reduced the size of the instances the most: 
		see Table \ref{table-redsize}. 
			\item It needs very little additional memory, precisely $O(nm).$ For details, and the 
		Matlab code, see Appendix \ref{section:matlabcode}.
			\item It is very accurate and stable: it is as accurate as Cholesky factorization, which works in machine precision.
		Sieve-SDP  is also easily implemented in a {\em safe mode}: see Subsection \ref{subsect:safemode}. 
			\item It is the simplest: the core Matlab code consists of only 65 lines.
	\end{itemize}

\noindent The code is available from 
	\begin{center}
		%\begin{verbatim}
		\url{https://github.com/unc-optimization/SieveSDP}
		%\end{verbatim}
	\end{center}

{\bf Acknowledgements} We thank the Technical Editor and the referees for their helpful comments.
The second author, G\'abor Pataki, 
is supported by the National Science Foundation, award DMS-1817272. 
The third author, Quoc Tran-Dinh, is supported in part by the National Science Foundation, award  DMS-1619884.
We are very grateful to Erling Andersen at Mosek for running several SDPs, and explaining the results; to Joachim Dahl at Mosek for helpful discussions on converting SDPs, and for providing his conversion code;  to Didier Henrion and Kim-Chuan Toh for providing us with some of the datasets; to Frank Permenter and Johan L$\ddot{o}$fberg  for helpful comments; to Oktay G\"unl\"uk for helping us to invent the name ``Sieve-SDP";
and to Hans Mittelmann for helping us with some of the large-scale SDPs.

\appendix

%%%%%%%
% \section{Very detailed results}
%%%%%%%

\section{Very detailed results}\label{section:verydetailed}

We now give very detailed computational results 
on all problems, separately for the five datasets. We only report on problems that were reduced by at least one of the five preprocessors.

In all tables the first column gives the number of the SDP, the second gives the name,  and the third gives the names of the 
preprocessing methods.

The next two columns describe the size of the problem. 
The entry ``f; l; s'' describes the size of the variables of the problem, where 
 \begin{itemize}
 	\item the number ``f'' is  the number of {\em free} variables;
 	\item the number ``l'' is the number of {\em linear nonnegative} variables;
 	\item the number, or numbers ``s'' describes the size of the {\em semidefinite}  variable blocks, possibly with multiplicity.
 \end{itemize}
For example, $3; 5; 6$ means that a problem has $3$ free variables; 
$5$ linear nonnegative variables; and a semidefinite matrix variable block  
of order $6.$ 
The tuple  $3; 5; 6, 5_3$ means that a problem has $3$ free variables; 
$5$ linear  nonnegative variables; and four semidefinite matrix variable blocks,  which are of order $6, 5, 5, 5, \,$ respectively.
The number $m$ is the number of constraints.

In the next three columns we put information about the preprocessors. In the column 
``red." we put $1$, if a preprocessor reduced a problem, and $0$ if it did not. In this column 
under Sieve-SDP we put the same entries, except when  Sieve-SDP actually proved infeasibility. In that case we entered  ``infeas" there. 
The number $\text{t}_{\rm prep}$ is the time spent on preprocessing;
the number $\text{t}_{\rm conv}$ is the time spent on converting 
from Mosek format to Sedumi format and back (for the methods pd1, pd2, dd1, and dd2).

In the next four columns we show how Mosek performed.
 In the column ``infeas" we have a $1$ if Mosek detected infeasibility, and $0$ if it did  not. 
The column 
obj (P, D) shows the objective values (primal and dual, respectively). 
The column DIMACS contains the \textit{largest} absolute value of the DIMACS errors.
The number $\text{t}_{\rm sol}$ is the time spent on solving the SDP. 

In the last column we show {\em help codes},  which show whether a preprocessor helped or hurt to solve an SDP.
  Although the help codes can be deduced from the previous columns, they still help to
 quickly evaluate the preprocessors. %how the codes performed. 
 A positive help code means that a preprocessor helped, and a negative one means that it hurt.
 
In detail, let us recall from 
 Subsection \ref{subsection-criteria-comparison} that 
 $\text{DIMACS}_\text{before}$ [$\text{DIMACS}_\text{after}$] is the absolute value of the DIMACS error 
 that  is largest in absolute value {before} [{after}] preprocessing.
We let  $\text{obj}_{\text{before}}$ and $\text{obj}_{\text{after}}$ be the primal objective values before and after preprocessing, respectively.

Given this notation, 
  	\begin{itemize}
 	 	\item the help code 	is $1, \,$ if
 	    \bit 
 		\item	 Sieve-SDP detects infeasibility, or 
 		\item Mosek does not detect infeasibility 
 		{\em before} preprocessing  but it does detect infeasibility {\em after} preprocessing; 
 	    \eit 
 	 	\item the help code is $-1, \,$ if 
 	    \bit 
 		\item Mosek detects infeasibility {\em before} preprocessing  but does not detect infeasibility 		{\em after} preprocessing; 
 	     \eit 
     \item  the help code is $2, \,$ if 
     \bit  
     	\item   it is not $\pm 1$ and preprocessing improved the DIMACS error, i.e., 
     	$$\text{DIMACS}_\text{before} > {10}^{-6} \, ~~{\rm and}~~ \, \dfrac{\text{DIMACS}_\text{after}}{\text{DIMACS}_\text{before}} < \dfrac{1}{10}; 
     	$$ 
     	\eit
     	\item the help code is $-2, \, $ if 
     	\bit 
     		\item  	it is not $\pm 1$ and preprocessing made the DIMACS error worse, i.e., 
     	$$\text{DIMACS}_{\text{after}} > {10}^{-6} \, ~~{\rm and}~~ \, \dfrac{\text{DIMACS}_\text{after}}{\text{DIMACS}_\text{before}} > 10; 
     	$$ 
     \eit  
      \end{itemize} 
  \begin{itemize}
  	\item the help code is  $3,$ if preprocessing shifted the objective value, i.e., 
  	\bit 
  		\item if help codes $\pm 1$ and $-2$ do not apply, and   
  		$$
  		\dfrac{|\text{obj}_{\text{before}} - \text{obj}_{\text{after}}|}{1 +  |\text{obj}_{\text{before}}|} > 10^{-6}; 
  		$$
  	\eit 
\end{itemize} 
\begin{itemize}
	\item the help code is MM if a code ran out of memory or crashed. 
\end{itemize}

%\begin{landscape}
% \subsection{Detailed results on the Permenter-Parrilo (PP) dataset}
%\newpage

\subsection{Detailed results on the Permenter-Parrilo (PP) dataset}

This dataset   has $68$ problems. From these $59$ problems were reduced by at least one of the five methods.
%
%\begin{footnotesize}
%\begin{table}[H]
%\begin{scriptsize}
%\begin{center}
%\setlength{\tabcolsep}{2pt}
%\begin{longtabu}[!htbp]{lllllllllllll}
\begin{table}[htbp!]
%\vspace{-5ex}
\setlength{\tabcolsep}{3pt}
\begin{center}
\resizebox{\textwidth}{!}{
\begin{tabular}{rllcrcrrccrrc}
%\begin{tabular}{| c | l | c | c c | c c c | c c c c | c | }
\hline\noalign{\smallskip}
No. & name & prep. method & f; l; s & m & red. & $\text{t}_\text{prep}$ & $\text{t}_\text{conv}$ & infeas & obj (P, D) & DIMACS & $\text{t}_\text{sol}$ & help\\
\noalign{\smallskip}\hline\noalign{\smallskip}
\multirowcell{6}{1} & \multirowcell{6}{CompactDim2R1} & none & $0;3;3$ & 5 & & & & 0 & 3.79e+06, 4.20e+06 & 2.22e+01 & 3.02 & \\ 
& & pd1 & $0;3;1$ & 3 & 1 & 0.05 & 0.00 & 1 & 0.00e+00, 1.00e+00 & 7.07e-01 & 0.64 & 1 \\ 
& & pd2 & $0;3;1$ & 3 & 1 & 0.04 & 0.00 & 1 & 0.00e+00, 1.00e+00 & 7.07e-01 & 0.69 & 1 \\ 
& & dd1 & & & 0 & 0.03 & 0.00 & & & & & \\ 
& & dd2 & & & 0 & 0.03 & 0.00 & & & & & \\ 
& & Sieve-SDP & & & infeas & 0.01 & & & & & 0.00 & 1 \\ 
\noalign{\smallskip}\hline\noalign{\smallskip}
\multirowcell{6}{2} & \multirowcell{6}{CompactDim2R2} & none & $0;0;6,{3}_{3}$ & 14 & & & & 0 & 6.41e-10, 6.81e-10 & 7.07e-01 & 3.16 & \\
\hline 
\end{tabular}}
\end{center}
\end{table}
\begin{table}[H]
\begin{center}
\setlength{\tabcolsep}{3pt}
\resizebox{\textwidth}{!}{	
\begin{tabular}{rllcrcrrccrrc}
\hline\noalign{\smallskip}
No. & name & prep. method & f; l; s & m & red. & $\text{t}_\text{prep}$ & $\text{t}_\text{conv}$ & infeas & obj (P, D) & DIMACS & $\text{t}_\text{sol}$ & help\\
\noalign{\smallskip}\hline\noalign{\smallskip}
\multirowcell{6}{5} & \multirowcell{6}{CompactDim2R5} & none & $0;0;21,{15}_{3}$ & 65 & & & & 0 & 1.50e+00, 1.50e+00 & 1.83e-07 & 2.05 & \\ 
& & pd1 & $0;0;{1}_{3}$ & 2 & 1 & 0.25 & 0.00 & 1 & 1.00e+00, 2.00e+00 & 7.07e-01 & 1.06 & 1 \\ 
& & pd2 & $0;0;{1}_{3}$ & 2 & 1 & 0.27 & 0.00 & 1 & 1.00e+00, 2.00e+00 & 7.07e-01 & 1.07 & 1 \\ 
& & dd1 & & & 0 & 0.02 & 0.00 & & & & & \\ 
& & dd2 & & & 0 & 0.05 & 0.00 & & & & & \\ 
& & Sieve-SDP & & & infeas & 0.03 & & & & & 0.00 & 1 \\ 
\noalign{\smallskip}\hline\noalign{\smallskip}
\multirowcell{6}{6} & \multirowcell{6}{CompactDim2R6} & none & $0;0;28,{21}_{3}$ & 90 & & & & 0 & 1.50e+00, 1.50e+00 & 2.70e-07 & 2.06 & \\
\hline 
\end{tabular}}
\end{center}
\end{table}
\begin{table}[H]
\begin{center}
\setlength{\tabcolsep}{3pt}
\resizebox{\textwidth}{!}{	
\begin{tabular}{rllcrcrrccrrc}
\hline\noalign{\smallskip}
No. & name & prep. method & f; l; s & m & red. & $\text{t}_\text{prep}$ & $\text{t}_\text{conv}$ & infeas & obj (P, D) & DIMACS & $\text{t}_\text{sol}$ & help\\
\noalign{\smallskip}\hline\noalign{\smallskip}
\multirowcell{6}{13} & \multirowcell{6}{Example3} & none & $0;0;3$ & 4 & & & & 0 & 3.33e-01, 3.33e-01 & 6.90e-02 & 1.79 & \\ 
& & pd1 & $0;0;2$ & 1 & 1 & 0.03 & 0.00 & 0 & 1.17e-07, 1.69e-07 & 5.14e-08 & 1.44 & 2,3 \\ 
& & pd2 & $0;0;2$ & 1 & 1 & 0.03 & 0.00 & 0 & 1.17e-07, 1.69e-07 & 5.14e-08 & 1.48 & 2,3 \\ 
& & dd1 & $3;0;2$ & 4 & 1 & 0.02 & 0.01 & 0 & 4.73e-15, 1.82e-14 & 2.75e-14 & 1.00 & 2,3 \\ 
& & dd2 & $3;0;2$ & 4 & 1 & 0.03 & 0.00 & 0 & 4.73e-15, 1.82e-14 & 2.75e-14 & 0.99 & 2,3 \\ 
& & Sieve-SDP & $0;0;2$ & 1 & 1 & 0.01 & & 0 & 1.17e-07, 1.69e-07 & 5.14e-08 & 1.59 & 2,3 \\ 
\noalign{\smallskip}\hline\noalign{\smallskip}
\multirowcell{6}{14} & \multirowcell{6}{Example4} & none & $0;0;3$ & 3 & & & & 1 & 0.00e+00, 3.74e-07 & 5.00e-01 & 1.43 & \\
\hline 
\end{tabular}}
\end{center}
\end{table}
\begin{table}[H]
\begin{center}
\setlength{\tabcolsep}{3pt}
\resizebox{\textwidth}{!}{	
\begin{tabular}{rllcrcrrccrrc}
\hline\noalign{\smallskip}
No. & name & prep. method & f; l; s & m & red. & $\text{t}_\text{prep}$ & $\text{t}_\text{conv}$ & infeas & obj (P, D) & DIMACS & $\text{t}_\text{sol}$ & help\\
\noalign{\smallskip}\hline\noalign{\smallskip}
\multirowcell{6}{21} & \multirowcell{6}{RandGen8} & none & $0;0;60$ & 40 & & & & 0 & 5.41e-09, 2.44e-09 & 9.31e-08 & 0.83 & \\ 
& & pd1 & & & 0 & 0.04 & 0.01 & & & & & \\ 
& & pd2 & & & 0 & 0.22 & 0.01 & & & & & \\ 
& & dd1 & & & 0 & 0.02 & 0.01 & & & & & \\ 
& & dd2 & $1269;0;33$ & 40 & 1 & 0.33 & 0.02 & 0 & 2.15e-15, 2.78e-19 & 6.90e-14 & 1.05 &  \\ 
& & Sieve-SDP & $0;0;30$ & 20 & 1 & 0.03 & & 0 & 1.52e-09, 6.33e-10 & 9.04e-09 & 0.69 &  \\ 
\noalign{\smallskip}\hline\noalign{\smallskip} 
\multirowcell{6}{22} & \multirowcell{6}{copos\_1} & none & $0;0;35$ & 210 & & & & 0 & 0.00e+00, 1.11e-08 & 4.40e-07 & 0.66 & \\
\hline 
\end{tabular}}
\end{center}
\end{table}
\begin{table}[H]
\begin{center}
\setlength{\tabcolsep}{3pt}
\resizebox{\textwidth}{!}{	
\begin{tabular}{rllcrcrrccrrc}
\hline\noalign{\smallskip}
No. & name & prep. method & f; l; s & m & red. & $\text{t}_\text{prep}$ & $\text{t}_\text{conv}$ & infeas & obj (P, D) & DIMACS & $\text{t}_\text{sol}$ & help\\
\noalign{\smallskip}\hline\noalign{\smallskip} 
\multirowcell{6}{29} & \multirowcell{6}{horn2} & none & $0;0;4$ & 7 & & & & 0 & 0.00e+00, 6.69e-13 & 9.06e-13 & 2.03 & \\ 
& & pd1 & & & 0 & 0.02 & 0.00 & & & & & \\ 
& & pd2 & $0;0;2$ & 3 & 1 & 0.06 & 0.00 & 0 & 0.00e+00, 0.00e+00 & 1.57e-16 & 0.99 &  \\ 
& & dd1 & & & 0 & 0.02 & 0.00 & & & & & \\ 
& & dd2 & & & 0 & 0.03 & 0.00 & & & & & \\ 
& & Sieve-SDP & & & 0 & 0.01 & & & & & & \\
\noalign{\smallskip}\hline\noalign{\smallskip}
\multirowcell{6}{30} & \multirowcell{6}{horn3} & none & $0;0;10$ & 28 & & & & 0 & 0.00e+00, 1.46e-07 & 8.62e-07 & 2.00 & \\
\hline 
\end{tabular}}
\end{center}
\end{table}
\begin{table}[H]
\begin{center}
\setlength{\tabcolsep}{3pt}
\resizebox{\textwidth}{!}{	
\begin{tabular}{rllcrcrrccrrc}
\hline\noalign{\smallskip}
No. & name & prep. method & f; l; s & m & red. & $\text{t}_\text{prep}$ & $\text{t}_\text{conv}$ & infeas & obj (P, D) & DIMACS & $\text{t}_\text{sol}$ & help\\
\noalign{\smallskip}\hline\noalign{\smallskip} 
\multirowcell{6}{38} & \multirowcell{6}{leverage\_limit} & none & $0;18100;{151}_{100},{30}_{100}$ & 68195 & & & & 0 & -8.75e+01, -8.75e+01 & 1.53e-05 & 278.60 & \\ 
& & pd1 & & & 0 & 2.10 & 0.17 & & & & & \\ 
& & pd2 & $0;18100;{151}_{99},121,{30}_{100}$ & 67700 & 1 & 120.98 & 7.87 & 0 & -8.75e+01, -8.75e+01 & 5.63e-06 & 150.78 & 3 \\ 
& & dd1 & $958500;18100;{61}_{100},{30}_{100}$ & 68195 & 1 & 3.87 & 7.20 & 0 & -8.75e+01, -8.75e+01 & 2.45e-05 & 250.27 &  \\ 
& & dd2 & $1193505;18100;{1}_{99},31$ & 68195 & 1 & 291.58 & 1.39 & -1 & -3.35e+00, 0.00e+00 & 1.03e+01 & 1.97 & -2 \\ 
& & Sieve-SDP & $0;18100;{143}_{97},{141}_{3},{26}_{98},{25}_{2}$ & 56196 & 1 & 253.43 & & 0 & -8.74e+01, -8.74e+01 & 1.73e-05 & 179.26 & 3 \\
\noalign{\smallskip}\hline\noalign{\smallskip} 
\multirowcell{6}{39} & \multirowcell{6}{long\_only} & none & $0;9000;{91}_{100},{30}_{100}$ & 59095 & & & & 0 & -4.13e+01, -4.13e+01 & 5.23e-06 & 373.38 & \\
\hline 
\end{tabular}}
\end{center}
\end{table}
\begin{table}[H]
\begin{center}
\setlength{\tabcolsep}{3pt}
\resizebox{\textwidth}{!}{	
\begin{tabular}{rllcrcrrccrrc}
\hline\noalign{\smallskip}
No. & name & prep. method & f; l; s & m & red. & $\text{t}_\text{prep}$ & $\text{t}_\text{conv}$ & infeas & obj (P, D) & DIMACS & $\text{t}_\text{sol}$ & help\\
\noalign{\smallskip}\hline\noalign{\smallskip} 
\multirowcell{6}{48} & \multirowcell{6}{unboundDim1R7} & none & $0;0;8,{7}_{2}$ & 14 & & & & 0 & -1.00e+00, -1.00e+00 & 5.11e-08 & 2.82 & \\ 
& & pd1 & $0;0;{1}_{2}$ & 1 & 1 & 0.21 & 0.00 & 0 & 0.00e+00, 0.00e+00 & 0.00e+00 & 0.45 & 3 \\ 
& & pd2 & $0;0;{1}_{2}$ & 1 & 1 & 0.23 & 0.00 & 0 & 0.00e+00, 0.00e+00 & 0.00e+00 & 0.44 & 3 \\ 
& & dd1 & & & 0 & 0.03 & 0.00 & & & & & \\ 
& & dd2 & & & 0 & 0.03 & 0.00 & & & & & \\ 
& & Sieve-SDP & $0;0;{1}_{2}$ & 1 & 1 & 0.02 & & 0 & 0.00e+00, 0.00e+00 & 0.00e+00 & 2.32 & 3 \\ 
\noalign{\smallskip}\hline\noalign{\smallskip} 
\multirowcell{6}{49} & \multirowcell{6}{unboundDim1R8} & none & $0;0;9,{8}_{2}$ & 16 & & & & 0 & -1.00e+00, -1.00e+00 & 5.43e-08 & 2.29 & \\
\end{tabular}}
\end{center}
\vspace{-5ex}
\end{table}
\begin{table}[H]
\begin{center}
\setlength{\tabcolsep}{3pt}
\resizebox{\textwidth}{!}{	
\begin{tabular}{rllcrcrrccrrc}
\hline\noalign{\smallskip}
No. & name & prep. method & f; l; s & m & red. & $\text{t}_\text{prep}$ & $\text{t}_\text{conv}$ & infeas & obj (P, D) & DIMACS & $\text{t}_\text{sol}$ & help\\
\noalign{\smallskip}\hline\noalign{\smallskip} 
\multirowcell{6}{56} & \multirowcell{6}{wei\_wagner\_W3\_PlusE} & none & $0;0;9$ & 38 & & & & 0 & 0.00e+00, -9.18e-09 & 5.53e-08 & 1.98 & \\ 
& & pd1 & & & 0 & 0.02 & 0.00 & & & & & \\ 
& & pd2 & $0;0;5$ & 15 & 1 & 0.06 & 0.00 & 0 & 0.00e+00, -7.21e-09 & 3.21e-08 & 1.02 &  \\ 
& & dd1 & & & 0 & 0.02 & 0.00 & & & & & \\ 
& & dd2 & & & 0 & 0.03 & 0.00 & & & & & \\ 
& & Sieve-SDP & & & 0 & 0.00 & & & & & & \\ 
\noalign{\smallskip}\hline\noalign{\smallskip} 
\multirowcell{6}{57} & \multirowcell{6}{wei\_wagner\_nP\_minus\_1\_24} & none & $0;0;12$ & 64 & & & & 0 & 0.00e+00, -5.50e-09 & 8.80e-08 & 2.03 & \\
\end{tabular}}
\end{center}
\vspace{-5ex}
\end{table}
%

%\end{center}
%\end{scriptsize}
%\end{table}

%\end{footnotesize}
% \subsection{Detailed results on the Mittelmann dataset}

\subsection{Detailed results on the Mittelmann  dataset}

This dataset  has $31$ problems. From these $8$ problems were reduced by at least one of the five methods. There were $5$ problems on which pd2 or dd2 ran out of memory or crashed.
%
%\begin{footnotesize}
%\begin{scriptsize}
%\begin{center}
%	\setlength{\tabcolsep}{2pt}
%%
\begin{table}[H]
%\vspace{-4ex}
\begin{center}
\setlength{\tabcolsep}{3pt}
\resizebox{\textwidth}{!}{	
\begin{tabular}{rllcrcrrccrrc}
%\begin{longtabu}[!htbp]{lllllllllllll}
\hline\noalign{\smallskip}
No. & name & prep. method & f; l; s & m & red. & $\text{t}_\text{prep}$ & $\text{t}_\text{conv}$ & infeas & obj (P, D) & DIMACS & $\text{t}_\text{sol}$ & help\\
\noalign{\smallskip}\hline\noalign{\smallskip}
\multirowcell{6}{1} & \multirowcell{6}{diamond\_patch} & none & $0;0;5477$ & 5478 & & & & 0 & 1.63e+01, 1.63e+01 & 3.56e-04 & 10854.97 & \\ 
& & pd1 & & & 0 & 31.05 & 0.06 & & & & & \\ 
& & pd2 & & & & MM & & & & & & MM \\ 
& & dd1 & & & 0 & 27.94 & 0.06 & & & & & \\ 
& & dd2 & & & 0 & 3008.86 & 0.06 & & & & & \\ 
& & Sieve-SDP & & & 0 & 1.12 & & & & & & \\ 
\noalign{\smallskip}\hline\noalign{\smallskip} 
\multirowcell{6}{2} & \multirowcell{6}{e\_moment\_stable\_17\_0.5\_2\_2} & none & $0;342;171,{18}_{17}$ & 5984 & & & & 0 & -1.98e-01, -1.98e-01 & 1.14e-05 & 38.53 & \\ 
& & pd1 & $0;342;{18}_{18}$ & 1139 & 1 & 0.71 & 0.16 & 0 & -1.98e-01, -1.98e-01 & 8.44e-06 & 1.64 &  \\ 
& & pd2 & $0;342;{18}_{18}$ & 1139 & 1 & 0.86 & 0.13 & 0 & -1.98e-01, -1.98e-01 & 8.44e-06 & 1.66 &  \\ 
& & dd1 & & & 0 & 0.08 & 0.01 & & & & & \\ 
& & dd2 & & & 0 & 0.34 & 0.01 & & & & & \\ 
& & Sieve-SDP & $0;342;{18}_{18}$ & 1139 & 1 & 0.52 & & 0 & -1.98e-01, -1.98e-01 & 8.44e-06 & 1.63 &  \\ 
\noalign{\smallskip}\hline\noalign{\smallskip} 
\multirowcell{6}{3} & \multirowcell{6}{ice\_2.0} & none & $0;0;8113$ & 8113 & & & & 0 & 6.81e+03, 6.81e+03 & 4.58e-07 & 17680.82 & \\ 
& & pd1 & & & 0 & 65.58 & 0.01 & & & & & \\ 
& & pd2 & & & & MM & & & & & & MM \\ 
& & dd1 & & & 0 & 64.03 & 0.01 & & & & & \\ 
& & dd2 & & & & MM & & & & & & MM \\ 
& & Sieve-SDP & & & 0 & 0.80 & & & & & & \\ 
\noalign{\smallskip}\hline 
\end{tabular}}
\end{center}
\end{table}
\begin{table}[H]
%\vspace{-4ex}
\begin{center}
\setlength{\tabcolsep}{3pt}
\resizebox{\textwidth}{!}{	
\begin{tabular}{rllcrcrrccrrc}
%\begin{longtabu}[!htbp]{lllllllllllll}
\hline\noalign{\smallskip}
No. & name & prep. method & f; l; s & m & red. & $\text{t}_\text{prep}$ & $\text{t}_\text{conv}$ & infeas & obj (P, D) & DIMACS & $\text{t}_\text{sol}$ & help\\
\noalign{\smallskip}\hline\noalign{\smallskip}
\multirowcell{6}{4} & \multirowcell{6}{G60\_mb} & none & $0;0;7000$ & 7001 & & & & 0 & 1.93e+03, 1.93e+03 & 6.64e-05 & 29138.79 & \\ 
& & pd1 & & & 0 & 107.66 & 10.87 & & & & & \\ 
& & pd2 & & & & MM & & & & & & MM \\ 
& & dd1 & & & 0 & 72.76 & 10.87 & & & & & \\ 
& & dd2 & & & & MM & & & & & & MM \\ 
& & Sieve-SDP & & & 0 & 22.42 & & & & & & \\ 
\noalign{\smallskip}\hline\noalign{\smallskip}
\multirowcell{6}{5} & \multirowcell{6}{maxG60} & none & $0;0;7000$ & 7000 & & & & 0 & -1.52e+04, -1.52e+04 & 6.73e-07 & 5217.88 & \\ 
& & pd1 & & & 0 & 47.47 & 0.01 & & & & & \\ 
& & pd2 & & & & MM & & & & & & MM \\ 
& & dd1 & & & 0 & 45.65 & 0.01 & & & & & \\ 
& & dd2 & & & & MM & & & & & & MM \\ 
& & Sieve-SDP & & & 0 & 0.47 & & & & & & \\ 
\noalign{\smallskip}\hline\noalign{\smallskip} 
\multirowcell{6}{6} & \multirowcell{6}{neu3} & none & $0;2;418$ & 7364 & & & & 0 & 7.10e-08, 1.12e-08 & 2.01e-06 & 153.03 & \\ 
& & pd1 & $0;2;87$ & 1152 & 1 & 0.94 & 0.11 & 0 & 4.69e-08, 3.50e-08 & 1.94e-07 & 3.01 & 2 \\ 
& & pd2 & $0;2;87$ & 1152 & 1 & 5.41 & 0.10 & 0 & 4.69e-08, 3.50e-08 & 1.94e-07 & 2.97 & 2 \\ 
& & dd1 & & & 0 & 0.16 & 0.02 & & & & & \\ 
& & dd2 & & & 0 & 2.29 & 0.02 & & & & & \\ 
& & Sieve-SDP & $0;2;87$ & 1152 & 1 & 2.34 & & 0 & 4.69e-08, 3.50e-08 & 1.94e-07 & 2.99 & 2 \\ 
\noalign{\smallskip}\hline\noalign{\smallskip} 
\multirowcell{6}{7} & \multirowcell{6}{neu3g} & none & $0;0;462$ & 8007 & & & & 0 & 4.58e-08, -2.89e-09 & 8.67e-07 & 151.22 & \\ 
& & pd1 & $0;0;87$ & 1151 & 1 & 1.32 & 0.11 & 0 & 8.91e-08, 5.65e-08 & 2.91e-07 & 3.00 &  \\ 
& & pd2 & $0;0;87$ & 1151 & 1 & 10.68 & 0.11 & 0 & 8.91e-08, 5.65e-08 & 2.91e-07 & 3.09 &  \\ 
& & dd1 & & & 0 & 0.19 & 0.03 & & & & & \\ 
& & dd2 & & & 0 & 2.66 & 0.03 & & & & & \\ 
& & Sieve-SDP & $0;0;87$ & 1151 & 1 & 2.26 & & 0 & 8.91e-08, 5.65e-08 & 2.91e-07 & 3.03 &  \\ 
\noalign{\smallskip}\hline\noalign{\smallskip} 
\multirowcell{6}{8} & \multirowcell{6}{p\_auss2\_3.0} & none & $0;0;9115$ & 9115 & & & & 0 & 8.62e+03, 8.62e+03 & 2.36e-07 & 25651.19 & \\ 
& & pd1 & & & 0 & 93.91 & 0.02 & & & & & \\ 
& & pd2 & & & & MM & & & & & & MM \\ 
& & dd1 & & & 0 & 97.11 & 0.02 & & & & & \\ 
& & dd2 & & & & MM & & & & & & MM \\ 
& & Sieve-SDP & & & 0 & 0.76 & & & & & & \\ 
\noalign{\smallskip}\hline\noalign{\smallskip} 
\multirowcell{6}{9} & \multirowcell{6}{rose13} & none & $0;0;105$ & 2379 & & & & 0 & 1.20e+01, 1.20e+01 & 1.65e-06 & 7.63 & \\ 
& & pd1 & $0;0;92$ & 1911 & 1 & 0.11 & 0.14 & 0 & 1.20e+01, 1.20e+01 & 4.86e-07 & 5.26 &  \\ 
& & pd2 & $0;0;80$ & 1523 & 1 & 0.51 & 0.11 & 0 & 1.20e+01, 1.20e+01 & 1.98e-07 & 2.94 &  \\ 
& & dd1 & & & 0 & 0.05 & 0.01 & & & & & \\ 
& & dd2 & & & 0 & 0.12 & 0.01 & & & & & \\ 
& & Sieve-SDP & $0;0;92$ & 1911 & 1 & 0.39 & & 0 & 1.20e+01, 1.20e+01 & 4.86e-07 & 5.28 &  \\ 
\noalign{\smallskip}\hline\noalign{\smallskip} 
\multirowcell{6}{10} & \multirowcell{6}{rose15} & none & $0;2;135$ & 3860 & & & & 0 & -3.11e-06, -2.94e-06 & 1.83e-05 & 19.47 & \\ 
& & pd1 & $0;2;121$ & 3181 & 1 & 0.08 & 0.24 & 0 & -3.52e-07, -1.52e-07 & 5.07e-05 & 11.73 & 3 \\ 
& & pd2 & $0;2;107$ & 2593 & 1 & 0.66 & 0.19 & 0 & -1.59e-09, -1.57e-09 & 1.10e-08 & 5.74 & 2,3 \\ 
& & dd1 & & & 0 & 0.07 & 0.00 & & & & & \\ 
& & dd2 & & & 0 & 0.18 & 0.00 & & & & & \\ 
& & Sieve-SDP & $0;2;121$ & 3181 & 1 & 0.52 & & 0 & -3.52e-07, -1.52e-07 & 5.07e-05 & 11.71 & 3 \\ 
\noalign{\smallskip}\hline\noalign{\smallskip}
\multirowcell{6}{11} & \multirowcell{6}{taha1a} & none & $0;0;252,{56}_{3},{126}_{10}$ & 3002 & & & & 0 & -1.00e+00, -1.00e+00 & 9.39e-07 & 37.54 & \\ 
& & pd1 & $0;0;126,{56}_{3},{126}_{10}$ & 2001 & 1 & 10.57 & 0.72 & 0 & -1.00e+00, -1.00e+00 & 1.20e-07 & 21.55 &  \\ 
& & pd2 & $0;0;126,{56}_{3},{126}_{10}$ & 2001 & 1 & 18.98 & 0.75 & 0 & -1.00e+00, -1.00e+00 & 1.20e-07 & 21.50 &  \\ 
& & dd1 & & & 0 & 0.21 & 0.06 & & & & & \\ 
& & dd2 & & & 0 & 21.47 & 0.06 & & & & & \\ 
& & Sieve-SDP & $0;0;126,{56}_{3},{126}_{10}$ & 2001 & 1 & 1.75 & & 0 & -1.00e+00, -1.00e+00 & 1.20e-07 & 21.70 &  \\ 
\noalign{\smallskip}\hline 
\end{tabular}}
\end{center}
\end{table}
\begin{table}[H]
\begin{center}
\setlength{\tabcolsep}{3pt}
\resizebox{\textwidth}{!}{	
\begin{tabular}{rllcrcrrccrrc}
\hline\noalign{\smallskip}
No. & name & prep. method & f; l; s & m & red. & $\text{t}_\text{prep}$ & $\text{t}_\text{conv}$ & infeas & obj (P, D) & DIMACS & $\text{t}_\text{sol}$ & help\\
\noalign{\smallskip}\hline\noalign{\smallskip} 
\multirowcell{6}{12} & \multirowcell{6}{taha1b} & none & $0;3;286,{66}_{20}$ & 8007 & & & & 0 & -7.73e-01, -7.73e-01 & 1.59e-07 & 148.99 & \\ 
& & pd1 & $0;3;{66}_{21}$ & 3002 & 1 & 13.97 & 0.87 & 0 & -7.73e-01, -7.73e-01 & 1.32e-07 & 34.29 &  \\ 
& & pd2 & $0;3;{66}_{21}$ & 3002 & 1 & 18.37 & 0.85 & 0 & -7.73e-01, -7.73e-01 & 1.32e-07 & 33.03 &  \\ 
& & dd1 & & & 0 & 0.16 & 0.04 & & & & & \\ 
& & dd2 & & & 0 & 1.82 & 0.04 & & & & & \\ 
& & Sieve-SDP & $0;3;{66}_{21}$ & 3002 & 1 & 1.97 & & 0 & -7.73e-01, -7.73e-01 & 1.32e-07 & 32.97 &  \\ 
\noalign{\smallskip}\hline\noalign{\smallskip} 
\multirowcell{6}{13} & \multirowcell{6}{taha1c} & none & $0;0;462,{126}_{3},{252}_{10}$ & 6187 & & & & 0 & -1.00e+00, -1.00e+00 & 3.12e-07 & 314.61 & \\ 
& & pd1 & $0;0;252,{126}_{3},{252}_{10}$ & 4367 & 1 & 148.36 & 2.11 & 0 & -1.00e+00, -1.00e+00 & 4.37e-07 & 178.22 &  \\ 
& & pd2 & $0;0;252,{126}_{3},{252}_{10}$ & 4367 & 1 & 187.99 & 2.01 & 0 & -1.00e+00, -1.00e+00 & 4.37e-07 & 177.80 &  \\ 
& & dd1 & & & 0 & 0.75 & 0.25 & & & & & \\ 
& & dd2 & & & 0 & 156.72 & 0.25 & & & & & \\ 
& & Sieve-SDP & $0;0;252,{126}_{3},{252}_{10}$ & 4367 & 1 & 10.85 & & 0 & -1.00e+00, -1.00e+00 & 4.37e-07 & 182.86 &  \\ 
%\end{longtabu}
\noalign{\smallskip}\hline 
\end{tabular}}
\vspace{-5ex}
\end{center}
\end{table}
%

%\end{center}
%\end{scriptsize}
%\end{footnotesize}

% \subsection{Detailed results on the Dressler-Illiman-de Wolff (DIW)  dataset}

\subsection{Detailed results on the Dressler-Illiman-de Wolff (DIW)  dataset}

This is a collection of 155 SDP relaxations from polynomial optimization generated by Gloptipoly 3 based on the paper \cite{dressler2018approach}. All problems were reduced by at least one preprocessor. 
% 
%\begin{footnotesize}
%\begin{scriptsize}
%\begin{center}
%	\setlength{\tabcolsep}{2pt}
%%
\begin{table}[H]
%\vspace{-4ex}
\begin{center}
\setlength{\tabcolsep}{3pt}
\resizebox{\textwidth}{!}{	
\begin{tabular}{rllcrcrrccrrc}
%\hline\noalign{\smallskip}
%No. & name & prep. method & f; l; s & m & red. & $\text{t}_\text{prep}$ & $\text{t}_\text{conv}$ & infeas & obj (P, D) & DIMACS & $\text{t}_\text{sol}$ & help\\
%
%\begin{longtabu}[!htbp]{lllllllllllll}
\hline\noalign{\smallskip}
No. & name & prep. method & f; l; s & m & red. & $\text{t}_\text{prep}$ & $\text{t}_\text{conv}$ & infeas & obj (P, D) & DIMACS & $\text{t}_\text{sol}$ & help\\
\noalign{\smallskip}\hline\noalign{\smallskip}
\multirowcell{6}{1} & \multirowcell{6}{ex3.3\_order4} & none & $0;1;15$ & 44 & & & & 0 & 3.54e-10, 3.56e-10 & 6.87e-01 & 1.69 & \\ 
& & pd1 & $0;1;2$ & 3 & 1 & 0.10 & 0.00 & 1 & 0.00e+00, 5.00e-01 & 5.27e-01 & 0.63 & 1 \\ 
& & pd2 & $0;1;2$ & 3 & 1 & 0.35 & 0.00 & 1 & 0.00e+00, 5.00e-01 & 5.27e-01 & 0.73 & 1 \\ 
& & dd1 & & & 0 & 0.04 & 0.00 & & & & & \\ 
& & dd2 & & & 0 & 0.04 & 0.00 & & & & & \\ 
& & Sieve-SDP & & & infeas & 0.03 & & & & & 0.00 & 1 \\ 
\noalign{\smallskip}\hline\noalign{\smallskip} 
\multirowcell{6}{2} & \multirowcell{6}{ex3.3\_order5} & none & $0;0;21,3$ & 65 & & & & 0 & 6.16e-02, 6.16e-02 & 1.18e-06 & 1.25 & \\
\hline 
\end{tabular}}
\end{center}
\end{table}
\begin{table}[H]
\begin{center}
\setlength{\tabcolsep}{3pt}
\resizebox{\textwidth}{!}{	
\begin{tabular}{rllcrcrrccrrc}
\hline\noalign{\smallskip}
No. & name & prep. method & f; l; s & m & red. & $\text{t}_\text{prep}$ & $\text{t}_\text{conv}$ & infeas & obj (P, D) & DIMACS & $\text{t}_\text{sol}$ & help\\
\noalign{\smallskip}\hline\noalign{\smallskip}
\multirowcell{6}{5} & \multirowcell{6}{ex3.3\_order8} & none & $0;0;45,15$ & 152 & & & & 0 & 6.16e-02, 6.16e-02 & 2.45e-06 & 1.46 & \\ 
& & pd1 & $0;0;12,6$ & 38 & 1 & 0.25 & 0.01 & 0 & 6.16e-02, 6.16e-02 & 5.23e-08 & 0.65 & 2,3 \\ 
& & pd2 & $0;0;12,6$ & 38 & 1 & 0.37 & 0.01 & 0 & 6.16e-02, 6.16e-02 & 5.23e-08 & 0.67 & 2,3 \\ 
& & dd1 & & & 0 & 0.02 & 0.00 & & & & & \\ 
& & dd2 & & & 0 & 0.04 & 0.00 & & & & & \\ 
& & Sieve-SDP & $0;0;12,6$ & 38 & 1 & 0.07 & & 0 & 6.16e-02, 6.16e-02 & 5.23e-08 & 0.58 & 2,3 \\ 
\noalign{\smallskip}\hline\noalign{\smallskip} 
\multirowcell{6}{6} & \multirowcell{6}{ex3.3\_order9} & none & $0;0;55,21$ & 189 & & & & 0 & 6.16e-02, 6.16e-02 & 4.63e-06 & 0.94 & \\
\hline 
\end{tabular}}
\end{center}
\end{table}
\begin{table}[H]
\begin{center}
\setlength{\tabcolsep}{3pt}
\resizebox{\textwidth}{!}{	
\begin{tabular}{rllcrcrrccrrc}
\hline\noalign{\smallskip}
No. & name & prep. method & f; l; s & m & red. & $\text{t}_\text{prep}$ & $\text{t}_\text{conv}$ & infeas & obj (P, D) & DIMACS & $\text{t}_\text{sol}$ & help\\
\noalign{\smallskip}\hline\noalign{\smallskip} 
\multirowcell{6}{12} & \multirowcell{6}{ex3.3\_order15} & none & $0;0;136,78$ & 495 & & & & 0 & 6.16e-02, 6.16e-02 & 1.06e-05 & 2.72 & \\ 
& & pd1 & $0;0;22,16$ & 75 & 1 & 0.73 & 0.01 & 0 & 6.16e-02, 6.16e-02 & 1.47e-07 & 1.06 & 2,3 \\ 
& & pd2 & $0;0;22,16$ & 75 & 1 & 2.89 & 0.02 & 0 & 6.16e-02, 6.16e-02 & 1.47e-07 & 1.06 & 2,3 \\ 
& & dd1 & & & 0 & 0.05 & 0.00 & & & & & \\ 
& & dd2 & & & 0 & 0.24 & 0.00 & & & & & \\ 
& & Sieve-SDP & $0;0;22,16$ & 75 & 1 & 0.47 & & 0 & 6.16e-02, 6.16e-02 & 1.47e-07 & 1.10 & 2,3 \\ 
\noalign{\smallskip}\hline\noalign{\smallskip} 
\multirowcell{6}{13} & \multirowcell{6}{ex3.3\_order16} & none & $0;0;153,91$ & 560 & & & & 0 & 6.16e-02, 6.16e-02 & 1.18e-05 & 3.64 & \\
\hline 
\end{tabular}}
\end{center}
\end{table}
\begin{table}[H]
\begin{center}
\setlength{\tabcolsep}{3pt}
\resizebox{\textwidth}{!}{	
\begin{tabular}{rllcrcrrccrrc}
\hline\noalign{\smallskip}
No. & name & prep. method & f; l; s & m & red. & $\text{t}_\text{prep}$ & $\text{t}_\text{conv}$ & infeas & obj (P, D) & DIMACS & $\text{t}_\text{sol}$ & help\\
\noalign{\smallskip}\hline\noalign{\smallskip} 
\multirowcell{6}{20} & \multirowcell{6}{ex4.1\_order5} & none & $0;0;21,6$ & 65 & & & & 0 & 3.08e-09, 3.12e-09 & 8.29e-01 & 0.94 & \\ 
& & pd1 & $0;0;1$ & 1 & 1 & 0.10 & 0.00 & 1 & 0.00e+00, 3.00e+00 & 7.50e-01 & 0.39 & 1 \\ 
& & pd2 & $0;0;1$ & 1 & 1 & 0.11 & 0.00 & 1 & 0.00e+00, 3.00e+00 & 7.50e-01 & 0.34 & 1 \\ 
& & dd1 & & & 0 & 0.02 & 0.00 & & & & & \\ 
& & dd2 & & & 0 & 0.03 & 0.00 & & & & & \\ 
& & Sieve-SDP & & & infeas & 0.02 & & & & & 0.00 & 1 \\ 
\noalign{\smallskip}\hline\noalign{\smallskip} 
\multirowcell{6}{21} & \multirowcell{6}{ex4.1\_order6} & none & $0;0;28,10$ & 90 & & & & 0 & 1.00e+00, 1.00e+00 & 6.31e-07 & 0.85 & \\
\hline 
\end{tabular}}
\end{center}
\end{table}
\begin{table}[H]
\begin{center}
\setlength{\tabcolsep}{3pt}
\resizebox{\textwidth}{!}{	
\begin{tabular}{rllcrcrrccrrc}
\hline\noalign{\smallskip}
No. & name & prep. method & f; l; s & m & red. & $\text{t}_\text{prep}$ & $\text{t}_\text{conv}$ & infeas & obj (P, D) & DIMACS & $\text{t}_\text{sol}$ & help\\
\noalign{\smallskip}\hline\noalign{\smallskip}
\multirowcell{6}{27} & \multirowcell{6}{ex4.1\_order12} & none & $0;0;91,55$ & 324 & & & & 0 & 1.00e+00, 1.00e+00 & 2.29e-06 & 1.06 & \\ 
& & pd1 & $0;0;1$ & 1 & 1 & 0.50 & 0.00 & 1 & 0.00e+00, 3.00e+00 & 7.50e-01 & 0.34 & 1 \\ 
& & pd2 & $0;0;1$ & 1 & 1 & 1.02 & 0.00 & 1 & 0.00e+00, 3.00e+00 & 7.50e-01 & 0.41 & 1 \\ 
& & dd1 & & & 0 & 0.02 & 0.00 & & & & & \\ 
& & dd2 & & & 0 & 0.11 & 0.00 & & & & & \\ 
& & Sieve-SDP & & & infeas & 0.16 & & & & & 0.00 & 1 \\ 
\noalign{\smallskip}\hline\noalign{\smallskip} 
\multirowcell{6}{28} & \multirowcell{6}{ex4.1\_order13} & none & $0;0;105,66$ & 377 & & & & 0 & 1.00e+00, 1.00e+00 & 6.83e-06 & 1.00 & \\
\hline 
\end{tabular}}
\end{center}
\end{table}
\begin{table}[H]
\begin{center}
\setlength{\tabcolsep}{3pt}
\resizebox{\textwidth}{!}{	
\begin{tabular}{rllcrcrrccrrc}
\hline\noalign{\smallskip}
No. & name & prep. method & f; l; s & m & red. & $\text{t}_\text{prep}$ & $\text{t}_\text{conv}$ & infeas & obj (P, D) & DIMACS & $\text{t}_\text{sol}$ & help\\
\noalign{\smallskip}\hline\noalign{\smallskip}
\multirowcell{6}{35} & \multirowcell{6}{ex4.1\_order20} & none & $0;0;231,171$ & 860 & & & & 0 & 1.00e+00, 1.00e+00 & 5.65e-06 & 8.86 & \\ 
& & pd1 & $0;0;1$ & 1 & 1 & 2.68 & 0.02 & 1 & 0.00e+00, 3.00e+00 & 7.50e-01 & 0.35 & 1 \\ 
& & pd2 & $0;0;1$ & 1 & 1 & 12.71 & 0.02 & 1 & 0.00e+00, 3.00e+00 & 7.50e-01 & 0.34 & 1 \\ 
& & dd1 & & & 0 & 0.07 & 0.02 & & & & & \\ 
& & dd2 & & & 0 & 0.87 & 0.02 & & & & & \\ 
& & Sieve-SDP & & & infeas & 0.88 & & & & & 0.00 & 1 \\ 
\noalign{\smallskip}\hline\noalign{\smallskip}
\multirowcell{6}{36} & \multirowcell{6}{ex4.2\_order4} & none & $0;1;15$ & 44 & & & & 0 & 1.00e-09, 1.01e-09 & 7.07e-01 & 0.81 & \\
\hline 
\end{tabular}}
\end{center}
\end{table}
\begin{table}[H]
\begin{center}
\setlength{\tabcolsep}{3pt}
\resizebox{\textwidth}{!}{	
\begin{tabular}{rllcrcrrccrrc}
\hline\noalign{\smallskip}
No. & name & prep. method & f; l; s & m & red. & $\text{t}_\text{prep}$ & $\text{t}_\text{conv}$ & infeas & obj (P, D) & DIMACS & $\text{t}_\text{sol}$ & help\\
\noalign{\smallskip}\hline\noalign{\smallskip} 
\multirowcell{6}{43} & \multirowcell{6}{ex4.2\_order11} & none & $0;0;78,36$ & 275 & & & & 0 & 5.53e-01, 5.53e-01 & 1.61e-06 & 1.66 & \\ 
& & pd1 & $0;0;10,4$ & 28 & 1 & 0.49 & 0.01 & 0 & 5.53e-01, 5.53e-01 & 3.90e-08 & 0.80 & 2,3 \\ 
& & pd2 & $0;0;10,4$ & 28 & 1 & 0.82 & 0.01 & 0 & 5.53e-01, 5.53e-01 & 3.90e-08 & 0.73 & 2,3 \\ 
& & dd1 & & & 0 & 0.03 & 0.00 & & & & & \\ 
& & dd2 & & & 0 & 0.10 & 0.00 & & & & & \\ 
& & Sieve-SDP & $0;0;10,4$ & 28 & 1 & 0.16 & & 0 & 5.53e-01, 5.53e-01 & 3.90e-08 & 0.70 & 2,3 \\ 
\noalign{\smallskip}\hline\noalign{\smallskip}
\multirowcell{6}{44} & \multirowcell{6}{ex4.2\_order12} & none & $0;0;91,45$ & 324 & & & & 0 & 5.53e-01, 5.53e-01 & 2.77e-06 & 1.23 & \\
\hline 
\end{tabular}}
\end{center}
\end{table}
\begin{table}[H]
\begin{center}
\setlength{\tabcolsep}{3pt}
\resizebox{\textwidth}{!}{	
\begin{tabular}{rllcrcrrccrrc}
\hline\noalign{\smallskip}
No. & name & prep. method & f; l; s & m & red. & $\text{t}_\text{prep}$ & $\text{t}_\text{conv}$ & infeas & obj (P, D) & DIMACS & $\text{t}_\text{sol}$ & help\\
\noalign{\smallskip}\hline\noalign{\smallskip} 
\multirowcell{6}{50} & \multirowcell{6}{ex4.2\_order18} & none & $0;0;190,120$ & 702 & & & & 0 & 5.52e-01, 5.52e-01 & 1.70e-04 & 9.65 & \\ 
& & pd1 & $0;0;10,4$ & 28 & 1 & 2.01 & 0.02 & 0 & 5.53e-01, 5.53e-01 & 3.90e-08 & 0.62 & 2,3 \\ 
& & pd2 & $0;0;10,4$ & 28 & 1 & 8.06 & 0.02 & 0 & 5.53e-01, 5.53e-01 & 3.90e-08 & 0.65 & 2,3 \\ 
& & dd1 & & & 0 & 0.06 & 0.01 & & & & & \\ 
& & dd2 & & & 0 & 0.65 & 0.01 & & & & & \\ 
& & Sieve-SDP & $0;0;10,4$ & 28 & 1 & 0.78 & & 0 & 5.53e-01, 5.53e-01 & 3.90e-08 & 0.62 & 2,3 \\ 
\noalign{\smallskip}\hline\noalign{\smallskip} 
\multirowcell{6}{51} & \multirowcell{6}{ex4.2\_order19} & none & $0;0;210,136$ & 779 & & & & 0 & 5.52e-01, 5.52e-01 & 9.26e-04 & 11.00 & \\
\hline 
\end{tabular}}
\end{center}
\end{table}
\begin{table}[H]
\begin{center}
\setlength{\tabcolsep}{3pt}
\resizebox{\textwidth}{!}{	
\begin{tabular}{rllcrcrrccrrc}
\hline\noalign{\smallskip}
No. & name & prep. method & f; l; s & m & red. & $\text{t}_\text{prep}$ & $\text{t}_\text{conv}$ & infeas & obj (P, D) & DIMACS & $\text{t}_\text{sol}$ & help\\
\noalign{\smallskip}\hline\noalign{\smallskip} 
\multirowcell{6}{58} & \multirowcell{6}{ex4.3\_order7} & none & $0;0;120,56$ & 679 & & & & 0 & 1.60e+01, 1.60e+01 & 5.26e-06 & 4.45 & \\ 
& & pd1 & $0;0;1$ & 1 & 1 & 0.30 & 0.01 & 1 & 0.00e+00, 8.00e+00 & 8.89e-01 & 1.23 & 1 \\ 
& & pd2 & $0;0;1$ & 1 & 1 & 0.83 & 0.01 & 1 & 0.00e+00, 8.00e+00 & 8.89e-01 & 1.19 & 1 \\ 
& & dd1 & & & 0 & 0.06 & 0.01 & & & & & \\ 
& & dd2 & & & 0 & 0.23 & 0.01 & & & & & \\ 
& & Sieve-SDP & & & infeas & 0.17 & & & & & 0.00 & 1 \\ 
\noalign{\smallskip}\hline\noalign{\smallskip} 
\multirowcell{6}{59} & \multirowcell{6}{ex4.3\_order8} & none & $0;0;165,84$ & 968 & & & & 0 & 1.60e+01, 1.60e+01 & 5.17e-06 & 9.40 & \\
\hline 
\end{tabular}}
\end{center}
\end{table}
\begin{table}[H]
\begin{center}
\setlength{\tabcolsep}{3pt}
\resizebox{\textwidth}{!}{	
\begin{tabular}{rllcrcrrccrrc}
\hline\noalign{\smallskip}
No. & name & prep. method & f; l; s & m & red. & $\text{t}_\text{prep}$ & $\text{t}_\text{conv}$ & infeas & obj (P, D) & DIMACS & $\text{t}_\text{sol}$ & help\\
\noalign{\smallskip}\hline\noalign{\smallskip}
\multirowcell{6}{66} & \multirowcell{6}{ex4.3\_order15} & none & $0;0;816,560$ & 5455 & & & & 0 & 8.94e+00, 8.94e+00 & 1.76e-01 & 2010.36 & \\ 
& & pd1 & $0;0;1$ & 1 & 1 & 33.32 & 0.41 & 1 & 0.00e+00, 8.00e+00 & 8.89e-01 & 0.50 & 1 \\ 
& & pd2 & $0;0;1$ & 1 & 1 & 192.22 & 0.41 & 1 & 0.00e+00, 8.00e+00 & 8.89e-01 & 0.70 & 1 \\ 
& & dd1 & & & 0 & 0.99 & 0.41 & & & & & \\ 
& & dd2 & & & 0 & 303.20 & 0.41 & & & & & \\ 
& & Sieve-SDP & & & infeas & 16.20 & & & & & 0.00 & 1 \\ 
\noalign{\smallskip}\hline\noalign{\smallskip} 
\multirowcell{6}{67} & \multirowcell{6}{ex4.3\_order16} & none & $0;0;969,680$ & 6544 & & & & 0 & 7.95e+00, 7.95e+00 & 2.23e-01 & 3158.88 & \\
\hline 
\end{tabular}}
\end{center}
\end{table}
\begin{table}[H]
\begin{center}
\setlength{\tabcolsep}{3pt}
\resizebox{\textwidth}{!}{	
\begin{tabular}{rllcrcrrccrrc}
\hline\noalign{\smallskip}
No. & name & prep. method & f; l; s & m & red. & $\text{t}_\text{prep}$ & $\text{t}_\text{conv}$ & infeas & obj (P, D) & DIMACS & $\text{t}_\text{sol}$ & help\\
\noalign{\smallskip}\hline\noalign{\smallskip}
\multirowcell{6}{74} & \multirowcell{6}{ex4.4\_order5} & none & $0;0;56,35$ & 285 & & & & 0 & 5.66e-08, 5.85e-08 & 8.66e-01 & 2.11 & \\ 
& & pd1 & $0;0;1$ & 1 & 1 & 0.14 & 0.00 & 1 & 0.00e+00, 1.00e+00 & 5.00e-01 & 0.60 & 1 \\ 
& & pd2 & $0;0;1$ & 1 & 1 & 0.25 & 0.00 & 1 & 0.00e+00, 1.00e+00 & 5.00e-01 & 0.59 & 1 \\ 
& & dd1 & & & 0 & 0.03 & 0.00 & & & & & \\ 
& & dd2 & & & 0 & 0.14 & 0.00 & & & & & \\ 
& & Sieve-SDP & & & infeas & 0.04 & & & & & 0.00 & 1 \\ 
\noalign{\smallskip}\hline\noalign{\smallskip} 
\multirowcell{6}{75} & \multirowcell{6}{ex4.4\_order6} & none & $0;0;84,56$ & 454 & & & & 0 & 1.49e-08, 1.50e-08 & 8.66e-01 & 3.11 & \\
\hline 
\end{tabular}}
\end{center}
\end{table}
\begin{table}[H]
\begin{center}
\setlength{\tabcolsep}{3pt}
\resizebox{\textwidth}{!}{	
\begin{tabular}{rllcrcrrccrrc}
\hline\noalign{\smallskip}
No. & name & prep. method & f; l; s & m & red. & $\text{t}_\text{prep}$ & $\text{t}_\text{conv}$ & infeas & obj (P, D) & DIMACS & $\text{t}_\text{sol}$ & help\\
\noalign{\smallskip}\hline\noalign{\smallskip} 
\multirowcell{6}{82} & \multirowcell{6}{ex4.4\_order13} & none & $0;0;560,455$ & 3653 & & & & 0 & 1.45e-08, 2.31e-09 & 1.10e-06 & 238.42 & \\ 
& & pd1 & $0;0;1$ & 1 & 1 & 9.81 & 0.19 & 1 & 0.00e+00, 1.00e+00 & 5.00e-01 & 1.18 & 1 \\ 
& & pd2 & $0;0;1$ & 1 & 1 & 61.06 & 0.19 & 1 & 0.00e+00, 1.00e+00 & 5.00e-01 & 1.02 & 1 \\ 
& & dd1 & & & 0 & 0.49 & 0.19 & & & & & \\ 
& & dd2 & & & 0 & 100.91 & 0.19 & & & & & \\ 
& & Sieve-SDP & & & infeas & 5.35 & & & & & 0.00 & 1 \\ 
\noalign{\smallskip}\hline\noalign{\smallskip}
\multirowcell{6}{83} & \multirowcell{6}{ex4.4\_order14} & none & $0;0;680,560$ & 4494 & & & & 0 & -1.06e-08, -3.67e-08 & 1.87e-06 & 455.22 & \\
\hline 
\end{tabular}}
\end{center}
\end{table}
\begin{table}[H]
\begin{center}
\setlength{\tabcolsep}{3pt}
\resizebox{\textwidth}{!}{	
\begin{tabular}{rllcrcrrccrrc}
\hline\noalign{\smallskip}
No. & name & prep. method & f; l; s & m & red. & $\text{t}_\text{prep}$ & $\text{t}_\text{conv}$ & infeas & obj (P, D) & DIMACS & $\text{t}_\text{sol}$ & help\\
\noalign{\smallskip}\hline\noalign{\smallskip}
\multirowcell{6}{90} & \multirowcell{6}{ex5.4\_order5} & none & $0;0;21$ & 65 & & & & 0 & 2.28e+00, 2.28e+00 & 5.89e-07 & 0.90 & \\ 
& & pd1 & $0;0;10$ & 31 & 1 & 0.07 & 0.00 & 0 & 2.28e+00, 2.28e+00 & 6.97e-08 & 0.86 &  \\ 
& & pd2 & $0;0;10$ & 31 & 1 & 0.09 & 0.00 & 0 & 2.28e+00, 2.28e+00 & 6.97e-08 & 0.78 &  \\ 
& & dd1 & & & 0 & 0.02 & 0.00 & & & & & \\ 
& & dd2 & & & 0 & 0.02 & 0.00 & & & & & \\ 
& & Sieve-SDP & $0;0;10$ & 31 & 1 & 0.02 & & 0 & 2.28e+00, 2.28e+00 & 6.97e-08 & 0.77 &  \\ 
\noalign{\smallskip}\hline\noalign{\smallskip} 
\multirowcell{6}{91} & \multirowcell{6}{ex5.4\_order6} & none & $0;0;28$ & 90 & & & & 0 & 2.28e+00, 2.28e+00 & 1.89e-06 & 0.95 & \\
\hline 
\end{tabular}}
\end{center}
\end{table}
\begin{table}[H]
\begin{center}
\setlength{\tabcolsep}{3pt}
\resizebox{\textwidth}{!}{	
\begin{tabular}{rllcrcrrccrrc}
\hline\noalign{\smallskip}
No. & name & prep. method & f; l; s & m & red. & $\text{t}_\text{prep}$ & $\text{t}_\text{conv}$ & infeas & obj (P, D) & DIMACS & $\text{t}_\text{sol}$ & help\\
\noalign{\smallskip}\hline\noalign{\smallskip}
\multirowcell{6}{97} & \multirowcell{6}{ex5.4\_order12} & none & $0;0;91$ & 324 & & & & 0 & 2.28e+00, 2.28e+00 & 1.09e-05 & 2.52 & \\ 
& & pd1 & $0;0;10$ & 31 & 1 & 0.29 & 0.01 & 0 & 2.28e+00, 2.28e+00 & 6.97e-08 & 1.72 & 2,3 \\ 
& & pd2 & $0;0;10$ & 31 & 1 & 0.54 & 0.01 & 0 & 2.28e+00, 2.28e+00 & 6.97e-08 & 0.82 & 2,3 \\ 
& & dd1 & & & 0 & 0.03 & 0.00 & & & & & \\ 
& & dd2 & & & 0 & 0.16 & 0.00 & & & & & \\ 
& & Sieve-SDP & $0;0;10$ & 31 & 1 & 0.20 & & 0 & 2.28e+00, 2.28e+00 & 6.97e-08 & 1.20 & 2,3 \\ 
\noalign{\smallskip}\hline\noalign{\smallskip} 
\multirowcell{6}{98} & \multirowcell{6}{ex5.4\_order13} & none & $0;0;105$ & 377 & & & & 0 & 2.28e+00, 2.28e+00 & 4.58e-06 & 2.22 & \\
\hline 
\end{tabular}}
\end{center}
\end{table}
\begin{table}[H]
\begin{center}
\setlength{\tabcolsep}{3pt}
\resizebox{\textwidth}{!}{	
\begin{tabular}{rllcrcrrccrrc}
\hline\noalign{\smallskip}
No. & name & prep. method & f; l; s & m & red. & $\text{t}_\text{prep}$ & $\text{t}_\text{conv}$ & infeas & obj (P, D) & DIMACS & $\text{t}_\text{sol}$ & help\\
\noalign{\smallskip}\hline\noalign{\smallskip} 
\multirowcell{6}{104} & \multirowcell{6}{ex5.4\_order19} & none & $0;0;210$ & 779 & & & & 0 & 2.23e+00, 2.23e+00 & 2.09e-02 & 8.09 & \\ 
& & pd1 & $0;0;10$ & 31 & 1 & 0.98 & 0.01 & 0 & 2.28e+00, 2.28e+00 & 6.97e-08 & 0.77 & 2,3 \\ 
& & pd2 & $0;0;10$ & 31 & 1 & 4.38 & 0.01 & 0 & 2.28e+00, 2.28e+00 & 6.97e-08 & 0.77 & 2,3 \\ 
& & dd1 & & & 0 & 0.04 & 0.00 & & & & & \\ 
& & dd2 & & & 0 & 0.61 & 0.00 & & & & & \\ 
& & Sieve-SDP & $0;0;10$ & 31 & 1 & 0.67 & & 0 & 2.28e+00, 2.28e+00 & 6.97e-08 & 1.66 & 2,3 \\ 
\noalign{\smallskip}\hline\noalign{\smallskip}
\multirowcell{6}{105} & \multirowcell{6}{ex5.4\_order20} & none & $0;0;231$ & 860 & & & & 0 & 2.14e+00, 2.15e+00 & 5.08e-02 & 10.64 & \\
\hline 
\end{tabular}}
\end{center}
\end{table}
\begin{table}[H]
\begin{center}
\setlength{\tabcolsep}{3pt}
\resizebox{\textwidth}{!}{	
\begin{tabular}{rllcrcrrccrrc}
\hline\noalign{\smallskip}
No. & name & prep. method & f; l; s & m & red. & $\text{t}_\text{prep}$ & $\text{t}_\text{conv}$ & infeas & obj (P, D) & DIMACS & $\text{t}_\text{sol}$ & help\\
\noalign{\smallskip}\hline\noalign{\smallskip} 
\multirowcell{6}{111} & \multirowcell{6}{ex5.5\_order9} & none & $0;0;55$ & 189 & & & & 0 & 1.62e-01, 1.62e-01 & 1.33e-06 & 0.92 & \\ 
& & pd1 & $0;0;7$ & 20 & 1 & 0.14 & 0.00 & 0 & 3.07e-01, 3.07e-01 & 1.08e-07 & 0.77 & 2,3 \\ 
& & pd2 & $0;0;7$ & 20 & 1 & 0.22 & 0.00 & 0 & 3.07e-01, 3.07e-01 & 1.08e-07 & 0.76 & 2,3 \\ 
& & dd1 & & & 0 & 0.02 & 0.00 & & & & & \\ 
& & dd2 & & & 0 & 0.04 & 0.00 & & & & & \\ 
& & Sieve-SDP & $0;0;7$ & 20 & 1 & 0.05 & & 0 & 3.07e-01, 3.07e-01 & 1.08e-07 & 0.74 & 2,3 \\ 
\noalign{\smallskip}\hline\noalign{\smallskip}
\multirowcell{6}{112} & \multirowcell{6}{ex5.5\_order10} & none & $0;0;66$ & 230 & & & & 0 & 1.62e-01, 1.62e-01 & 1.21e-05 & 0.94 & \\
\hline 
\end{tabular}}
\end{center}
\end{table}
\begin{table}[H]
\begin{center}
\setlength{\tabcolsep}{3pt}
\resizebox{\textwidth}{!}{	
\begin{tabular}{rllcrcrrccrrc}
\hline\noalign{\smallskip}
No. & name & prep. method & f; l; s & m & red. & $\text{t}_\text{prep}$ & $\text{t}_\text{conv}$ & infeas & obj (P, D) & DIMACS & $\text{t}_\text{sol}$ & help\\
\noalign{\smallskip}\hline\noalign{\smallskip}
\multirowcell{6}{118} & \multirowcell{6}{ex5.5\_order16} & none & $0;0;153$ & 560 & & & & 0 & 1.62e-01, 1.62e-01 & 4.38e-05 & 2.43 & \\ 
& & pd1 & $0;0;7$ & 20 & 1 & 0.52 & 0.01 & 0 & 3.07e-01, 3.07e-01 & 1.08e-07 & 0.76 & 2,3 \\ 
& & pd2 & $0;0;7$ & 20 & 1 & 1.93 & 0.01 & 0 & 3.07e-01, 3.07e-01 & 1.08e-07 & 0.76 & 2,3 \\ 
& & dd1 & & & 0 & 0.03 & 0.01 & & & & & \\ 
& & dd2 & & & 0 & 0.22 & 0.01 & & & & & \\ 
& & Sieve-SDP & $0;0;7$ & 20 & 1 & 0.22 & & 0 & 3.07e-01, 3.07e-01 & 1.08e-07 & 0.77 & 2,3 \\ 
\noalign{\smallskip}\hline\noalign{\smallskip} 
\multirowcell{6}{119} & \multirowcell{6}{ex5.5\_order17} & none & $0;0;171$ & 629 & & & & 0 & 1.62e-01, 1.62e-01 & 6.83e-06 & 3.24 & \\
\hline 
\end{tabular}}
\end{center}
\end{table}
\begin{table}[H]
\begin{center}
\setlength{\tabcolsep}{3pt}
\resizebox{\textwidth}{!}{	
\begin{tabular}{rllcrcrrccrrc}
\hline\noalign{\smallskip}
No. & name & prep. method & f; l; s & m & red. & $\text{t}_\text{prep}$ & $\text{t}_\text{conv}$ & infeas & obj (P, D) & DIMACS & $\text{t}_\text{sol}$ & help\\
\noalign{\smallskip}\hline\noalign{\smallskip} 
\multirowcell{6}{125} & \multirowcell{6}{ex5.6\_order6} & none & $0;0;28$ & 90 & & & & 0 & 3.04e-01, 3.04e-01 & 3.45e-07 & 0.89 & \\ 
& & pd1 & $0;0;8$ & 21 & 1 & 0.07 & 0.00 & 0 & 3.04e-01, 3.04e-01 & 1.01e-08 & 0.78 &  \\ 
& & pd2 & $0;0;8$ & 21 & 1 & 0.09 & 0.00 & 0 & 3.04e-01, 3.04e-01 & 1.01e-08 & 0.86 &  \\ 
& & dd1 & & & 0 & 0.03 & 0.00 & & & & & \\ 
& & dd2 & & & 0 & 0.07 & 0.00 & & & & & \\ 
& & Sieve-SDP & $0;0;8$ & 21 & 1 & 0.05 & & 0 & 3.04e-01, 3.04e-01 & 1.01e-08 & 2.02 &  \\ 
\noalign{\smallskip}\hline\noalign{\smallskip}
\multirowcell{6}{126} & \multirowcell{6}{ex5.6\_order7} & none & $0;0;36$ & 119 & & & & 0 & 3.04e-01, 3.04e-01 & 3.21e-07 & 1.17 & \\
\hline 
\end{tabular}}
\end{center}
\end{table}
\begin{table}[H]
\begin{center}
\setlength{\tabcolsep}{3pt}
\resizebox{\textwidth}{!}{	
\begin{tabular}{rllcrcrrccrrc}
\hline\noalign{\smallskip}
No. & name & prep. method & f; l; s & m & red. & $\text{t}_\text{prep}$ & $\text{t}_\text{conv}$ & infeas & obj (P, D) & DIMACS & $\text{t}_\text{sol}$ & help\\
\noalign{\smallskip}\hline\noalign{\smallskip} 
\multirowcell{6}{132} & \multirowcell{6}{ex5.6\_order13} & none & $0;0;105$ & 377 & & & & 0 & 3.04e-01, 3.04e-01 & 7.94e-07 & 1.20 & \\ 
& & pd1 & $0;0;8$ & 21 & 1 & 0.28 & 0.00 & 0 & 3.04e-01, 3.04e-01 & 1.01e-08 & 0.77 &  \\ 
& & pd2 & $0;0;8$ & 21 & 1 & 0.71 & 0.00 & 0 & 3.04e-01, 3.04e-01 & 1.01e-08 & 0.84 &  \\ 
& & dd1 & & & 0 & 0.03 & 0.00 & & & & & \\ 
& & dd2 & & & 0 & 0.11 & 0.00 & & & & & \\ 
& & Sieve-SDP & $0;0;8$ & 21 & 1 & 0.13 & & 0 & 3.04e-01, 3.04e-01 & 1.01e-08 & 0.79 &  \\ 
\noalign{\smallskip}\hline\noalign{\smallskip} 
\multirowcell{6}{133} & \multirowcell{6}{ex5.6\_order14} & none & $0;0;120$ & 434 & & & & 0 & 3.04e-01, 3.04e-01 & 1.19e-06 & 1.31 & \\
\hline 
\end{tabular}}
\end{center}
\end{table}
\begin{table}[H]
\begin{center}
\setlength{\tabcolsep}{3pt}
\resizebox{\textwidth}{!}{	
\begin{tabular}{rllcrcrrccrrc}
\hline\noalign{\smallskip}
No. & name & prep. method & f; l; s & m & red. & $\text{t}_\text{prep}$ & $\text{t}_\text{conv}$ & infeas & obj (P, D) & DIMACS & $\text{t}_\text{sol}$ & help\\
\noalign{\smallskip}\hline\noalign{\smallskip}
\multirowcell{6}{139} & \multirowcell{6}{ex5.6\_order20} & none & $0;0;231$ & 860 & & & & 0 & 3.04e-01, 3.04e-01 & 3.35e-06 & 5.55 & \\ 
& & pd1 & $0;0;8$ & 21 & 1 & 1.48 & 0.01 & 0 & 3.04e-01, 3.04e-01 & 1.01e-08 & 1.02 & 2 \\ 
& & pd2 & $0;0;8$ & 21 & 1 & 6.57 & 0.01 & 0 & 3.04e-01, 3.04e-01 & 1.01e-08 & 0.91 & 2 \\ 
& & dd1 & & & 0 & 0.07 & 0.01 & & & & & \\ 
& & dd2 & & & 0 & 0.58 & 0.01 & & & & & \\ 
& & Sieve-SDP & $0;0;8$ & 21 & 1 & 0.53 & & 0 & 3.04e-01, 3.04e-01 & 1.01e-08 & 1.02 & 2 \\ 
\noalign{\smallskip}\hline\noalign{\smallskip}
\multirowcell{6}{140} & \multirowcell{6}{ex5.7\_order5} & none & $0;1;21$ & 65 & & & & 0 & 8.52e-09, -5.99e-09 & 7.10e-08 & 0.93 & \\
\hline 
\end{tabular}}
\end{center}
\end{table}
\begin{table}[H]
\begin{center}
\setlength{\tabcolsep}{3pt}
\resizebox{\textwidth}{!}{	
\begin{tabular}{rllcrcrrccrrc}
\hline\noalign{\smallskip}
No. & name & prep. method & f; l; s & m & red. & $\text{t}_\text{prep}$ & $\text{t}_\text{conv}$ & infeas & obj (P, D) & DIMACS & $\text{t}_\text{sol}$ & help\\
\noalign{\smallskip}\hline\noalign{\smallskip} 
\multirowcell{6}{146} & \multirowcell{6}{ex5.7\_order11} & none & $0;0;78,28$ & 275 & & & & 0 & 1.71e-08, -1.50e-08 & 2.11e-07 & 1.11 & \\ 
& & pd1 & $0;0;14,6$ & 45 & 1 & 0.36 & 0.01 & 0 & 1.73e-08, -3.62e-09 & 8.09e-08 & 0.91 &  \\ 
& & pd2 & $0;0;14,6$ & 45 & 1 & 0.75 & 0.01 & 0 & 1.73e-08, -3.62e-09 & 8.09e-08 & 0.93 &  \\ 
& & dd1 & & & 0 & 0.02 & 0.00 & & & & & \\ 
& & dd2 & & & 0 & 0.07 & 0.00 & & & & & \\ 
& & Sieve-SDP & $0;0;14,6$ & 45 & 1 & 0.17 & & 0 & 1.73e-08, -3.62e-09 & 8.09e-08 & 0.95 &  \\ 
\noalign{\smallskip}\hline\noalign{\smallskip} 
\multirowcell{6}{147} & \multirowcell{6}{ex5.7\_order12} & none & $0;0;91,36$ & 324 & & & & 0 & 5.85e-09, -6.21e-09 & 7.87e-08 & 1.04 & \\
\hline 
\end{tabular}}
\end{center}
\end{table}
\begin{table}[H]
\begin{center}
\setlength{\tabcolsep}{3pt}
\resizebox{\textwidth}{!}{	
\begin{tabular}{rllcrcrrccrrc}
\hline\noalign{\smallskip}
No. & name & prep. method & f; l; s & m & red. & $\text{t}_\text{prep}$ & $\text{t}_\text{conv}$ & infeas & obj (P, D) & DIMACS & $\text{t}_\text{sol}$ & help\\
\noalign{\smallskip}\hline\noalign{\smallskip} 
\multirowcell{6}{153} & \multirowcell{6}{ex5.7\_order18} & none & $0;0;190,105$ & 702 & & & & 0 & 1.76e-08, -1.75e-08 & 4.12e-07 & 3.01 & \\ 
& & pd1 & $0;0;14,6$ & 45 & 1 & 1.76 & 0.02 & 0 & 1.73e-08, -3.62e-09 & 8.09e-08 & 0.84 &  \\ 
& & pd2 & $0;0;14,6$ & 45 & 1 & 7.46 & 0.02 & 0 & 1.73e-08, -3.62e-09 & 8.09e-08 & 1.15 &  \\ 
& & dd1 & & & 0 & 0.07 & 0.01 & & & & & \\ 
& & dd2 & & & 0 & 0.58 & 0.01 & & & & & \\ 
& & Sieve-SDP & $0;0;14,6$ & 45 & 1 & 0.94 & & 0 & 1.73e-08, -3.62e-09 & 8.09e-08 & 1.36 &  \\ 
\noalign{\smallskip}\hline\noalign{\smallskip} 
\multirowcell{6}{154} & \multirowcell{6}{ex5.7\_order19} & none & $0;0;210,120$ & 779 & & & & 0 & 1.68e-08, -1.67e-08 & 4.25e-07 & 4.46 & \\ 
& & pd1 & $0;0;14,6$ & 45 & 1 & 2.49 & 0.02 & 0 & 1.73e-08, -3.62e-09 & 8.09e-08 & 1.09 &  \\ 
& & pd2 & $0;0;14,6$ & 45 & 1 & 10.47 & 0.02 & 0 & 1.73e-08, -3.62e-09 & 8.09e-08 & 0.89 &  \\ 
& & dd1 & & & 0 & 0.05 & 0.02 & & & & & \\ 
& & dd2 & & & 0 & 0.84 & 0.02 & & & & & \\ 
& & Sieve-SDP & $0;0;14,6$ & 45 & 1 & 1.11 & & 0 & 1.73e-08, -3.62e-09 & 8.09e-08 & 1.01 &  \\ 
\noalign{\smallskip}\hline\noalign{\smallskip} 
\multirowcell{6}{155} & \multirowcell{6}{ex5.7\_order20} & none & $0;0;231,136$ & 860 & & & & 0 & 1.54e-08, -1.69e-08 & 4.40e-07 & 5.84 & \\ 
& & pd1 & $0;0;14,6$ & 45 & 1 & 3.28 & 0.03 & 0 & 1.73e-08, -3.62e-09 & 8.09e-08 & 1.02 &  \\ 
& & pd2 & $0;0;14,6$ & 45 & 1 & 13.67 & 0.03 & 0 & 1.73e-08, -3.62e-09 & 8.09e-08 & 1.17 &  \\ 
& & dd1 & & & 0 & 0.08 & 0.02 & & & & & \\ 
& & dd2 & & & 0 & 1.04 & 0.02 & & & & & \\ 
& & Sieve-SDP & $0;0;14,6$ & 45 & 1 & 1.71 & & 0 & 1.73e-08, -3.62e-09 & 8.09e-08 & 0.86 &  \\ 
\noalign{\smallskip}\hline
%\end{longtabu}
\end{tabular}}
\end{center}
\vspace{-5ex}
\end{table}

%\end{center}
%\end{scriptsize}
%\end{footnotesize}

% \subsection{Detailed results on the Henrion-Toh dataset}

\subsection{Detailed results on the Henrion-Toh dataset}

This dataset has $98$ problems. From these $18$ problems were reduced by at least one of the five methods.
%
%\begin{footnotesize}
%\begin{scriptsize}
%\begin{center}
%\setlength{\tabcolsep}{2pt}
%%
\begin{table}[H]
%\vspace{-4ex}
\begin{center}
\setlength{\tabcolsep}{3pt}
\resizebox{\textwidth}{!}{	
\begin{tabular}{rllcrcrrccrrc}
%\begin{longtabu}[!htbp]{lllllllllllll}
\hline\noalign{\smallskip}
No. & name & prep. method & f; l; s & m & red. & $\text{t}_\text{prep}$ & $\text{t}_\text{conv}$ & infeas & obj (P, D) & DIMACS & $\text{t}_\text{sol}$ & help\\
\noalign{\smallskip}\hline\noalign{\smallskip}
\multirowcell{6}{1} & \multirowcell{6}{sedumi-brown} & none & $925;0;56$ & 461 & & & & 0 & -7.34e-09, 0.00e+00 & 3.75e-07 & 0.93 & \\ 
& & pd1 & $925;0;21$ & 251 & 1 & 0.15 & 0.05 & 0 & -9.33e-11, 0.00e+00 & 6.25e-09 & 0.80 &  \\ 
& & pd2 & $925;0;21$ & 251 & 1 & 0.23 & 0.02 & 0 & -9.33e-11, 0.00e+00 & 6.25e-09 & 0.78 &  \\ 
& & dd1 & & & 0 & 0.03 & 0.00 & & & & & \\ 
& & dd2 & & & 0 & 0.06 & 0.00 & & & & & \\ 
& & Sieve-SDP & $925;0;21$ & 251 & 1 & 0.04 & & 0 & -9.33e-11, 0.00e+00 & 6.25e-09 & 0.79 &  \\ 
\noalign{\smallskip}\hline\noalign{\smallskip} 
\multirowcell{6}{2} & \multirowcell{6}{sedumi-conform3} & none & $630;0;56$ & 285 & & & & 0 & 2.05e-08, 0.00e+00 & 4.54e-07 & 0.66 & \\ 
& & pd1 & $630;0;53$ & 273 & 1 & 0.04 & 0.02 & 0 & 2.51e-08, 0.00e+00 & 4.90e-07 & 0.71 &  \\ 
& & pd2 & $630;0;53$ & 273 & 1 & 0.11 & 0.03 & 0 & 2.51e-08, 0.00e+00 & 4.90e-07 & 0.73 &  \\ 
& & dd1 & & & 0 & 0.03 & 0.00 & & & & & \\ 
& & dd2 & & & 0 & 0.06 & 0.00 & & & & & \\ 
& & Sieve-SDP & $630;0;53$ & 273 & 1 & 0.02 & & 0 & 2.51e-08, 0.00e+00 & 4.90e-07 & 0.68 &  \\ 
\noalign{\smallskip}\hline\noalign{\smallskip} 
\multirowcell{6}{3} & \multirowcell{6}{sedumi-conform4} & none & $1890;0;84$ & 454 & & & & 0 & -2.51e-08, 0.00e+00 & 5.57e-06 & 0.84 & \\ 
& & pd1 & $1890;0;81$ & 442 & 1 & 0.07 & 0.04 & 0 & -6.47e-09, 0.00e+00 & 1.74e-06 & 0.81 &  \\ 
& & pd2 & $1890;0;81$ & 442 & 1 & 0.26 & 0.04 & 0 & -6.47e-09, 0.00e+00 & 1.74e-06 & 0.88 &  \\ 
& & dd1 & & & 0 & 0.02 & 0.00 & & & & & \\ 
& & dd2 & & & 0 & 0.08 & 0.00 & & & & & \\ 
& & Sieve-SDP & $1890;0;81$ & 442 & 1 & 0.02 & & 0 & -6.47e-09, 0.00e+00 & 1.74e-06 & 0.77 &  \\ 
\noalign{\smallskip}\hline 
\end{tabular}}
\end{center}
\end{table}
\begin{table}[H]
\begin{center}
\setlength{\tabcolsep}{3pt}
\resizebox{\textwidth}{!}{	
\begin{tabular}{rllcrcrrccrrc}
\hline\noalign{\smallskip}
No. & name & prep. method & f; l; s & m & red. & $\text{t}_\text{prep}$ & $\text{t}_\text{conv}$ & infeas & obj (P, D) & DIMACS & $\text{t}_\text{sol}$ & help\\
\noalign{\smallskip}\hline\noalign{\smallskip}
\multirowcell{6}{4} & \multirowcell{6}{sedumi-fp23} & none & $0;0;28,{7}_{13}$ & 209 & & & & 0 & 2.13e+02, 2.13e+02 & 3.97e-06 & 1.50 & \\ 
& & pd1 & $0;0;{7}_{14}$ & 83 & 1 & 0.09 & 0.02 & 0 & 2.13e+02, 2.13e+02 & 9.96e-07 & 1.60 & 3 \\ 
& & pd2 & $0;0;{7}_{14}$ & 83 & 1 & 0.13 & 0.02 & 0 & 2.13e+02, 2.13e+02 & 9.96e-07 & 1.44 & 3 \\ 
& & dd1 & & & 0 & 0.05 & 0.00 & & & & & \\ 
& & dd2 & & & 0 & 0.06 & 0.00 & & & & & \\ 
& & Sieve-SDP & $0;0;{7}_{14}$ & 83 & 1 & 0.04 & & 0 & 2.13e+02, 2.13e+02 & 9.96e-07 & 1.33 & 3 \\ 
\noalign{\smallskip}\hline\noalign{\smallskip} 
\multirowcell{6}{5} & \multirowcell{6}{sedumi-fp24} & none & $0;0;105,{14}_{35}$ & 2379 & & & & 0 & 1.95e+02, 1.95e+02 & 9.68e-08 & 6.36 & \\
\hline 
\end{tabular}}
\end{center}
\end{table}
\begin{table}[H]
\begin{center}
\setlength{\tabcolsep}{3pt}
\resizebox{\textwidth}{!}{	
\begin{tabular}{rllcrcrrccrrc}
\hline\noalign{\smallskip}
No. & name & prep. method & f; l; s & m & red. & $\text{t}_\text{prep}$ & $\text{t}_\text{conv}$ & infeas & obj (P, D) & DIMACS & $\text{t}_\text{sol}$ & help\\
\noalign{\smallskip}\hline\noalign{\smallskip} 
\multirowcell{6}{12} & \multirowcell{6}{sedumi-fp35} & none & $0;0;35,{20}_{8}$ & 164 & & & & 0 & 4.00e+00, 4.00e+00 & 5.76e-06 & 0.86 & \\ 
& & pd1 & $0;0;{20}_{8},10$ & 119 & 1 & 0.16 & 0.04 & 0 & 4.00e+00, 4.00e+00 & 5.66e-07 & 0.80 & 2,3 \\ 
& & pd2 & $0;0;{20}_{8},10$ & 119 & 1 & 0.36 & 0.05 & 0 & 4.00e+00, 4.00e+00 & 5.66e-07 & 0.85 & 2,3 \\ 
& & dd1 & & & 0 & 0.03 & 0.00 & & & & & \\ 
& & dd2 & & & 0 & 0.28 & 0.00 & & & & & \\ 
& & Sieve-SDP & $0;0;{20}_{8},10$ & 119 & 1 & 0.05 & & 0 & 4.00e+00, 4.00e+00 & 5.66e-07 & 0.89 & 2,3 \\ 
\noalign{\smallskip}\hline\noalign{\smallskip}
\multirowcell{6}{13} & \multirowcell{6}{sedumi-fp44} & none & $0;0;4,{3}_{2}$ & 6 & & & & 0 & 4.44e+02, 4.44e+02 & 4.67e-08 & 1.02 & \\
\end{tabular}}
\vspace{-3ex}
\end{center}
\end{table}

%\end{center}
%\end{scriptsize}
%\end{footnotesize}
%\end{landscape}

%%%%%%%
% \section{Core Matlab code}
%%%%%%%

\section{Core Matlab code}\label{section:matlabcode}

In this section we provide our core Matlab code of Sieve-SDP (not including input, output,  and dual solution recovery) 
with some comments.
In our code we physically delete rows and columns of the $A_i$ and of $C$ only at the very end. During the execution 
of the algorithm 
we only {\em mark}  
such rows, columns and constraints as deleted. % using a marker array.

We use two arrays to keep track of what has been marked deleted:
\begin{enumerate}
	\item The $m$-vector \texttt{undeleted}, whose $i$th entry is $1$ if 
	constraint $i$ {\em has not} been deleted, and $0$ if it {\em has} been deleted.
	\item The sparse array $\texttt{I} \in \{0, 1 \}^{n \times (m+1)}$ with entries defined as follows.
	\begin{enumerate}  
		\item For all $i$ and for $1 \leq j \leq m,$  
		$$
		I(i,j) \, = \, \left\{  \begin{array}{rll} 0, & {\rm if \, in} \, A_j \, {\rm the \,} i{\rm th \, row \, and \, column \, are \, all \, zero \, or \, have \, been \, deleted};  \\ 
		1,  & {\rm otherwise.} \end{array} 
		\right. \, 
		$$          
		\item For all $i,$ 
		$$
		I(i,m+1) \, = \, \left\{ \begin{array}{rll} 0, & {\rm if \, in} \, {\rm all} \, A_j \, {\rm the \,} i{\rm th \, row \, and \, column \,  have \, been \, deleted}; \\ 
		1, & {\rm otherwise}.  \end{array} 
		\right. \, 
		$$          
	\end{enumerate}
\end{enumerate}

\lstset{language=Matlab,%
	%basicstyle=\color{red},
	breaklines=true,%
	morekeywords={matlab2tikz},
	keywordstyle=\color{blue},%
	morekeywords=[2]{1}, keywordstyle=[2]{\color{blue}},
	identifierstyle=\color{black},%
	stringstyle=\color{mylilas},
	commentstyle=\color{mygreen},%
	showstringspaces=false,%without this there will be a symbol in the places where there is a space
	%numbers=left,%
	%numberstyle={\tiny \color{black}},% size of the numbers
	%numbersep=9pt, % this defines how far the numbers are from the text
	emph=[1]{for,end,break},emphstyle=[1]\color{blue}, %some words to emphasise
	%emph=[2]{word1,word2}, emphstyle=[2]{style},    
}
\lstset{basicstyle=\ttfamily\footnotesize,breaklines=true}

\begin{lstlisting}
function[Ar, br, Cr, info] = SieveSDP(A, b, C, EPS)
    
    % Inputs:
    %  A: n-by-n*m sparse matrix,
    %     which is m symmetric n-by-n matrices side by side
    %  b: the vector of rhs in R^m, and b <= 0;
    %  C: the objective coefficient n-by-n matrix;
    %  EPS: accuracy for safe mode, with default value eps
    % Outputs: 
    %   Ar, br, cr: Reduced data after preprocessing
    %   info: A structure containing preprocessing info
    
    if nargin < 4, EPS = eps; end
    sqrtEPS = sqrt(EPS);

    Ar = []; br = []; Cr = [];
    n = size(C, 1); m = length(b);
    I = true(n, m + 1); % initial nonzero indices
    for i = 1:m,
        I(:, i) = any(A(:, (n*(i - 1) + 1):(n*i)), 2);
    end

    not_done    = 1;          % 1 means preprocessing not done
    undeleted   = ones(m, 1); % keep track of deleted constraints
    constr_ind  = (1:m);      % indices or undeleted constraints
    mr          = m;          % reduced number of constraints
    info.infeas = 0;          % infeasibility detected?
    info.red    = 0;          % any reduction?
    
    bn = -sqrtEPS*max(1, norm(b, inf));
                     % b < 0 if b < -sqrt(epsilon)*max{1, ||b||}
    bz = bn*sqrtEPS; % b = 0 if -epsilon*max{1, ||b||} < b <= 0

    % Preprocessing
    while not_done
        not_done = 0;
        for ii = 1:mr
            i  = constr_ind(ii);
            Ii = I(:, i); % indicates undeleted vars in matrix i
            Ai = A(Ii, n*(i - 1) + find(Ii)); % nonzero submatrix
            Iaux = any(Ai, 2);
            if find(Iaux == false, 1),
                I(Ii, i) = Iaux; Ii = I(:, i); Ai = Ai(Iaux, Iaux);
            end
            if isempty(Ai)                    
                if b(i) < bn, info.infeas = 1; return; end
                    % Ai=0 and bi<0 => infeasible
                if b(i) > bz, undeleted(i) = 0; continue; end
                    % Ai=0 and bi=0 => reduce
            end
            if b(i) < bn
                [~, pd_check] = chol(Ai);
                if pd_check == 0, info.infeas = 1; return; end
                    % Ai pd and bi<0 => infeasible
            else
                if b(i) > bz
                    [~, pd_check] = chol(Ai);
                    if pd_check == 0
                        % Ai pd and bi=0 => reduce
                        I(Ii, :) = false; undeleted(i) = 0;
                        not_done = 1;
                    else
                        [~, nd_check] = chol(-Ai);
                        if nd_check == 0
                            % Ai nd and bi=0 => reduce
                            I(Ii, :) = false; undeleted(i) = 0;
                            not_done = 1;
                        end
                    end
                end
            end
        end
        constr_ind = find(undeleted); mr = length(constr_ind);
    end

    % Undeleted rows/columns are marked in I(:, m + 1)
    % Now do physical deletion
    if mr == m
        Ar = A; br = b; Cr = C; info.red = 0; return;
    end
    info.red = 1;
    I_nonzero = I(:, m + 1); nr = nnz(I_nonzero);
    Ar        = sparse(nr, nr*mr);
    for ii = 1:mr
        i = constr_ind(ii);
        Ar(:, (nr*(ii - 1) + 1):(nr*ii)) ...
            = A(I_nonzero, n*(i - 1) + find(I_nonzero));
    end
    br = b(constr_ind);
    Cr = C(I_nonzero, I_nonzero);
    
end
\end{lstlisting}

%%%%%%%
% \section{The DIMACS errors}
%%%%%%%

\section{The DIMACS errors}\label{section:dimacs}

For the sake of completeness in this section 
we describe the DIMACS errors, which are commonly used to measure the accuracy of  approximate solutions 
$X$ of \eqref{p} and  $y$ of (\ref{d}).

Define the operator ${\cal A}: \rad{m} \rightarrow \sym{n}$ and its adjoint as 
\begin{eqnarray}
\A(X) & = & (A_1 \bullet X, \dots, A_m  \bullet X), \nonumber\\
\A^*(y)  & = & \sum_{i=1}^m y_i A_i.  \nonumber
\end{eqnarray}
Suppose we are given an approximate solution $X$ of \eqref{p} and an approximate  solution 
$y$ of (\ref{d}). For brevity, define $Z = C - \A^*(y).$

Then the DIMACS error measures are defined as follows: 
\begin{eqnarray}
\err_1 & = & \dfrac{\norm{\A(X) - b}_2}{1 + \norm{b}_\infty}, \nonumber\\
\err_2 & = & \max \biggl\{   0, \dfrac{-\lambda_{\min}(X)}{ 1 + \norm{b}_\infty}  \biggr\}, \nonumber\\ 
\err_3 & = & \dfrac{\norm{ \A^*(y) - C  - Z}_F}{  1 + \norm{C}_\infty}, \nonumber\\
\err_4 & = & \max \biggl\{   0,  \dfrac{-\lambda_{\min}(Z)}{1 + \norm{C}_\infty} \biggr\},  \nonumber\\
\err_5 & = & \dfrac{b^{\top} y - C \bullet X}{1 + | C \bullet X |  + | b^{\top} y |}, \nonumber\\
\err_6 & = & \dfrac{Z \bullet X}{  1 + | C \bullet X| + | b^{\top} y| }. \nonumber
\end{eqnarray}
In the above equations we use the following notation. If $M = (m_{ij})  \in \sym{n},$ then we 
write $\norm{M}_F$ for the  Frobenius norm of $M$
and $\norm{M}_\infty$ for the infinity norm of $M, $ i.e., 
\begin{eqnarray*}
	\norm{M}_F & = & \sqrt{\sum_{i,j} m_{ij}^2} \\
	\norm{M}_\infty & = & \max_{i,j} | m_{ij} |. 
\end{eqnarray*}
We also write $\lambda_{\min}(M)$ for the smallest eigenvalue of $M.$

%%%%%%%
% \section{Dual solution recovery}
%%%%%%%

\section{Dual solution recovery}\label{section:dual_recovery}

In this section we address the following question: suppose we preprocessed the problem 
\eqref{p} by Sieve-SDP, then computed an optimal solution of the preprocessed SDP, $(P_\text{pre})$, and of its dual, $(D_\text{pre})$.
Can we compute an optimal solution  of the original primal \eqref{p} and of its dual \eqref{d}? 
	The answer to the first question (primal solution recovery) 
is easy, while the issue of dual solution recovery is much more subtle.

First let us look at primal solution recovery.
		  Since Sieve-SDP deletes rows and columns from the variable matrix $X$ 
		  that  are always zero anyway, if $X^{\rm pre}$
		  is an optimal solution of 
		$(P_{\rm pre}),$ 
		then by simply padding $X^{\rm pre}$ with zeroes we obtain an optimal solution of \eqref{p}.

Next we discuss dual solution recovery. For simplicity %To illustrate the main ideas, 
we first assume that Sieve-SDP performed just one iteration. 
		Further, we also assume that in the Basic Step (in 
		 Figure \ref{figure-basicstep})  it eliminated the constraint $A_1 \bullet X = 0, \, $ where 
	\begin{equation*}
	A_1 = \begin{pmatrix}
	D	& 0\\
	0	& 0
	\end{pmatrix}, 
	\end{equation*}
	with  $D \succ 0$ and we let $r$ be the order of $D.$ %assume that $D$ is of order, say, $r.$ %it is, say, of order $r.$

Next, let us  write out $(D_{\rm pre}):$ 
	\begin{equation*}\tag{$D_{\rm pre}$} 
	\left.\begin{array}{ll}
	\displaystyle\sup_{y}	& \displaystyle\sum_{i = 2}^m b_i y_i\\
	\text{s.t.}\	& C - \displaystyle\sum_{i = 2}^m y_i A_i \in \begin{pmatrix}
	\times	& \times\\
	\times	& \oplus
	\end{pmatrix},  
	\end{array}\right.
	\end{equation*}
	where the notation means that the lower right  $(n-r) \times (n-r)$ principal  
	 block of $C - \sum_{i = 2}^m y_i A_i $ is positive semidefinite, and the rest is  arbitrary. Thus clearly 
	\beq \label{ineq-sieve-2} 
	\val \eqref{d}  \leq \val(D_{\rm pre}), % \, {\rm holds},  %\leq \val(P_{\rm pre}) =  \val\eqref{p} \,\, 
	\eeq  
since $(D_\text{pre})$ has a  feasible region which is at least as large as that of (\ref{d}) (and usually it is larger).
Assume that $y^{\rm pre} = (y^{\rm pre}_2, \dots, y^{\rm pre}_m)$ 
	is an optimal solution of $(D_{\rm pre}).$
Our recovery procedure, which we call Basic-Recovery, fixes $y^{\rm pre}$ and 
seeks $y_1$ such that  $(y_1, y^{\rm pre})$ is feasible in \eqref{d}, i.e., 
	\begin{equation} \label{eqn-linesearch} 
y_1 A_1 +  \sum_{i = 2}^m y_i^\text{pre}  A_i \preceq C.  
	\end{equation} 
  We do this by a very basic linesearch: we first try the values $y_1 = 0, -1, \, $ and $-2.$ If these all fail, then we try $y_1 = -100.$ If we fail with $y_1 = -100, $ we stop; otherwise   we test $y_1 = -3, -4, \dots$ and find 
  the largest $y_1$ such that (\ref{eqn-linesearch}) holds.

To test whether \eqref{eqn-linesearch} holds, again Cholesky factorization comes into play: using it, we test whether  
$$ 
C - (y_1 A_1 +  \sum_{i = 2}^m y_i^\text{pre}  A_i) + 10^{-6}I \succ 0 \,\, {\rm holds.} 
$$ 

%To test whether \eqref{eqn-linesearch} holds,  Cholesky factorization again comes into play: we test the relation 
 % $$ 
 % C - (y_1 A_1 +  \sum_{i = 2}^m y_i^\text{pre}  A_i) + 10^{-6}I \succ 0. \,
  %$$ 
   
Basic-Recovery is inspired by the dual solution recovery procedure in \cite{PerPar:14}, 
  which builds on the ideas in \cite{pataki2017bad},   and it assumes that the 
  dual problem \eqref{d} is reduced\footnote{See Remark 
  	\ref{remark:PP-primal-vs-dual} about how the primal and dual are defined  in \cite{PerPar:14}.}. 

The procedure Basic-Recovery may fail. To see why, first assume it succeeds, i.e., it computes 
	a feasible solution of \eqref{d}. Since $y_1$ has zero objective coefficient in \eqref{d}, 
	this solution has objective value  
 $\val(D_\text{pre}),$ hence by inequality \eqref{ineq-sieve-2} it is optimal in \eqref{d}, thus 
 $\val \eqref{d}  = \val( D_\text{pre} ).$ 
 Conversely, if 
$\val \eqref{d} < \val(D_\text{pre}),  \, $ then Basic-Recovery {\em must}  fail.

\bex (Example \ref{ex2} continued) When we apply Sieve-SDP to the SDP     
(\ref{posgap}),  it deletes the first row and first column in all matrices and it also deletes the first constraint.

Let us write out $(D_{\rm pre})$ again for this problem (i.e., repeat the SDP \eqref{posgap-dual-pre-1}):   
  \beq \label{posgap-dual-pre} 
  \left.\ba{rl}
  \displaystyle\sup_{y_2}  &  y_2   \\
  \mathrm{s.t.} &  y_2 \bpx  1 & 0 \\  0 & 0 \epx \preceq   \bpx   1 & 0 \\  0 & 0 \epx,  
  \ena\right.
  \eeq
  whose optimal solution is $y^{\rm pre}_2 = 1.$

Thus, Basic-Recovery seeks $y_1$ such that 
$$
  \ba{rl}
  y_1 \bpx 1 & 0 & 0 \\ 0 & 0 & 0 \\ 0 & 0 & 0 \epx + \bpx 0 & 0 & 1 \\ 0 & 1 & 0 \\ 1 & 0 & 0 \epx \preceq 
  \bpx 1 & 0 & 0 \\ 0 & 1 & 0 \\ 0 & 0 & 0 \epx,
  \ena
  $$
  and clearly there is no such $y_1.$
\eex

One can construct more sophisticated examples in which $\val(D_\text{pre}) = \val \eqref{d}, \, $ 
but Basic-Recovery still fails.

We next look at dual solution recovery when Sieve-SDP deleted several constraints: then we run Basic-Recovery to find the corresponding $y_i$ sequentially. 
For simplicity 
assume that Sieve-SDP deleted constraints $1, 2, \dots, k$ and we found an optimal primal and dual 
solution of the resulting SDP (by Mosek). 
We then attempt  to find an optimal dual solution of the SDP obtained by deleting only constraints 
$1, \dots, k-1;$ then of the SDP obtained by deleting only constraints 
$1, \dots, k-2;$ and so on.

To conclude this section we make the point that dual solution recovery is much more difficult in SDP than in 
	linear programming. We thus implemented an ``ideal" recovery 
	 procedure, which we call Ideal-Recovery. It works   as follows.
Suppose  
	 $y^{\rm pre} =   (y^{\rm pre}_{k+1}, \dots,   y^{\rm pre}_{m})$ is an optimal dual 
	solution of the SDP obtained by deleting constraints $1, \dots, k.$  Ideal-Recovery  fixes $y^{\rm pre}, \, $ then calls Mosek to find a feasible solution $(y_1, \dots, y_k)$ of
		\begin{equation} \label{eqn-Ideal-recovery} 
	\sum_{i=1}^k y_i A_i + 	\sum_{i=k+1}^m y^{\rm pre}_i A_i  \preceq C. 
	\end{equation}
	Table \ref{tbl:solution_recovery} shows on how many instances pd1, pd2, Sieve-SDP+Basic-Recovery and Sieve-SDP+Ideal-Recovery succeeded. 
	(Note that they succeeded on overlapping, but 
	different problem sets, as a preprocessor may reduce an  SDP, while another preprocessor may not  reduce the same 
	SDP. We do not report results with dd1 and dd2, since they  reduced only very few instances.)
	
What do we mean by ``success"? 
 For  pd1 and pd2 it  means that their dual solution recovery code reported success. 
 For Sieve-SDP+Basic-Recovery it means that 
 it suceeded in every iteration: it computed the $y_i$ for every deleted constraint. 
  For Sieve-SDP+Ideal-Recovery  it means that Mosek did {\em not} report that 
  (\ref{eqn-Ideal-recovery}) is infeasible.

\begin{table}[!htbp]
%\vspace{-2ex}
\begin{center}
	\caption{Dual solution recovery by four methods}\label{tbl:solution_recovery}
	\setlength{\tabcolsep}{2pt}
	\resizebox{\textwidth}{!}{	
	\begin{tabular}{ l  r r r r r }
		\hline\noalign{\smallskip}
		Method & \# Reduced feasible  & \# Success & \# Failure & Success rate & Time (s) \\
		\noalign{\smallskip}\hline\noalign{\smallskip}
		pd1				& 137	& 23	& 114	& 16.8\%	& 154.75	\\
		pd2				& 158	& 39	& 119	& 24.7\%	& 172.13	\\
		Sieve-SDP + Basic-Recovery	& 143	& 25	& 118	& 17.5\%	& 12.62		\\
		Sieve-SDP + Ideal-Recovery  & 143	& 103	& 40	& 72.0\%	& 1313.57	\\
		\noalign{\smallskip}\hline
	\end{tabular}}
\end{center}
\vspace{-3ex}	
\end{table}

Next we made the criterion of ``success" more rigorous:  we redefined ``success" as returning a pair of primal-dual optimal solutions whose 
largest DIMACS error in absolute value is at most $10^{-6}.$ 
Table \ref{tbl:solution_recovery2} shows the results: 
now  Sieve-SDP+Basic-Recovery is the winner, as it beats the supposedly perfect Sieve-SDP+Ideal-Recovery 
 procedure.
 
\begin{table}[!htbp]
%\vspace{-2ex}
\begin{center}
	\caption{Dual solution recovery assuming the tightest standard for ``success"}\label{tbl:solution_recovery2}
	\setlength{\tabcolsep}{2pt}
	\resizebox{\textwidth}{!}{		
	\begin{tabular}{ l  r r r r r }
		\hline\noalign{\smallskip}
		Method & \# Reduced feasible & \# Success & \# Failure & Success rate & Time (s) \\
		\noalign{\smallskip}\hline\noalign{\smallskip}
		pd1				& 137	& 19	& 118	& 13.9\%	& 154.75	\\
		pd2				& 158	& 34	& 124	& 21.5\%	& 172.13	\\
		Sieve-SDP + Basic-Recovery		& 143	& 25	& 118	& 17.5\%	& 12.62		\\
		Sieve-SDP + Ideal-Recovery & 143	& 17	& 126	& 11.9\%	& 1313.57 \\
		\noalign{\smallskip}\hline
	\end{tabular}}
\end{center}
\vspace{-3ex}	
\end{table}

Nevertheless, none of the methods do very well, and Table \ref{tbl:solution_recovery2} shows that 
dual solution recovery in facial reduction remains a challenge, and an interesting area for further research.

%\bibliography{mpcrefs}   % name your BibTeX data base

%\begin{thebibliography}{}
% etc
%\end{thebibliography}

\end{document}